\documentclass[12pt]{amsart}
\usepackage{fontenc,amsmath,mathrsfs,amsfonts,amsthm, color,amssymb, MnSymbol, graphicx, hyperref, geroENG, multicol, enumerate, enumitem}
\usepackage[english]{babel}
\usepackage{times}[doublespace]

\usepackage[top=2.5cm, bottom=2.5cm, left=2.75cm, right=2.75cm]{geometry}
\usepackage{parskip}
\begin{document}

\title[Transcendence and normality via HCFs]{Transcendence and normality of complex numbers via Hurwitz continued fractions}
\author[Garc{\'i}a-Ramos]{Felipe García-Ramos}
\address{Felipe Garc{\'i}a-Ramos, Physics Institute, Universidad Aut{\'o}noma de San Luis Potos{\'i}, M{\'e}xico.
\newline
Faculty of Mathematics and Computer Science, Uniwersytet Jagielloński, Poland}
\email{fgramos@conahcyt.mx}

\author[Gonz{\'a}lez Robert]{Gerardo Gonz{\'a}lez Robert}
\address{Gerardo Gonz{\'a}lez Robert,  Department of Mathematical and Physical Sciences,  La Trobe University, Bendigo 3552, Australia. }
\email{G.Robert@latrobe.edu.au}

\author[Hussain]{Mumtaz Hussain}
\address{Mumtaz Hussain,  Department of Mathematical and Physical Sciences,  La Trobe University, Bendigo 3552, Australia. }
\email{m.hussain@latrobe.edu.au}


\begin{abstract}
We study the topological, dynamical, and descriptive set-theoretic properties of Hurwitz continued fractions. Hurwitz continued fractions associate an infinite sequence of Gaussian integers to every complex number which is not a Gaussian rational. The resulting space of sequences of Gaussian integers $\Omega$ is not closed. Using an iterative procedure, we show that $\Omega$ contains a natural subset whose closure $\overline{\mathsf{R}}$ encodes continued fraction expansions of complex numbers which are not Gaussian rationals. We prove that $(\overline{\mathsf{R}}, \sigma)$ is a subshift with a feeble specification property. As an application, we determine the rank in the Borel hierarchy of the set of Hurwitz normal numbers with respect to the complex Gauss measure. We also construct a family of complex transcendental numbers with bounded partial quotients.
\end{abstract}
 \maketitle
\setcounter{tocdepth}{1}

\section{Introduction}
A regular continued fraction is a finite or infinite expression of the form 
\begin{equation}\label{Eq:CF}
a_0 + \cfrac{1}{a_1 + \cfrac{1}{a_2 + \cfrac{1}{\ddots}}}    
\end{equation}
where all the terms $a_n$ are integers and $a_n\geq 1$ for any $n\geq 1$. Regular continued fractions are a powerful tool in number theory. They provide a number system where the expansion of a real number is infinite if and only if it is irrational, and eventually periodic if and only if the number is a quadratic irrational; for example,
\[
\sqrt2=1+\cfrac{1}{2+\cfrac{1}{2+\ddots}}.
\] 
From a Diophantine approximation perspective, regular continued fractions systematically give the best rational approximations to a real number (see, for example, the classic reference \cite{Khi1997}). We can study regular continued fractions from a dynamical perspective with the aid of the Gauss map $T_{\RE}:[0,1)\to [0,1)$ given by $T_{\RE}(0)=0$ and $T_{\RE}(x)=x^{-1}-\lfloor x^{-1}\rfloor$ for $x\neq 0$ where $\lfloor \cdot \rfloor$ is the floor function. This dynamical approach not only provides a natural framework to study regular continued fractions but also establishes connections with other fields such as harmonic analysis, probability theory, and descriptive set theory \cite{AirJacKwiMan2020, DiamondVaaler1986, Kaufman1980}.

There have been several attempts to obtain multidimensional continued fraction algorithms (see, for example, \cite{Karpenkov2022, Schweiger2000}). Still, the higher-dimensional theory, in comparison, is hardly as comprehensive as the one-dimensional theory. One such example is the continued fraction theory for the complex plane we focus on. In his 1887 paper \cite{Hur1887}, Hurwitz introduced a complex continued fraction algorithm where all the partial quotients $\sanu$ are Gaussian integers.  
Despite its age, the Hurwitz algorithm had not received extensive attention until recently, when interesting discoveries were made regarding its Diophantine approximation and dynamical properties \cite{BugGeroHus2023, Dani2015, DanNog2014, EiItoNak2019, HeXio2021-01}. Hurwitz showed that his continued fractions associate to each number $z\in \Cx\setminus \QU(i)$ a unique sequence $(a_n)_{n\geq 0}$ in $\Za[i]$ satisfying
\[
z
=
[a_0;a_1,a_2,\ldots]_{\Cx}
:=
\lim_{n\to\infty} a_0 + \cfrac{1}{a_1 + \cfrac{1}{\ddots + \cfrac{1}{a_n}}}.
\]
Like regular continued fractions, Hurwitz continued fractions can also be approached dynamically. The square 
\[
\mfF
:=
\left\{ z\in \Cx: -\frac{1}{2}\leq \real(z) < \frac{1}{2}, \; -\frac{1}{2}\leq \imag(z) < \frac{1}{2}\right\}
\]
plays the role of the unit interval, and the definition of the Hurwitz map $T$ is similar to that of $T_{\RE}$ (see Section \ref{SECTION:HCF}). We say a sequence $\sanu$ in $\Za[i]$ is \textit{valid} when $z=[0;a_1, a_2, \ldots]_{\Cx}$ for some $z\in \mfF\setminus \QU(i)$ and we denote by $\Omega$ the space of valid sequences. A central difference between regular and Hurwitz continued fractions is that, while $\Na^{\Na}$ is the space of valid sequences for regular continued fractions, the set $\Omega$ is not even closed in $\Za[i]^{\Na}$ (see Section \ref{SECTION:PFTHM1}). Working with closed sets is often crucial for developing a natural dynamical perspective. One solution might involve taking the closure of $\Omega$, but this approach can result in several intricate points that prove challenging to handle. Nevertheless, it turns out that there is a naturally closed subset related to $\Omega$ that still encapsulates its fundamental characteristics. To be more precise, for each $m\in\Na$ and $\bfb=(b_1,\ldots, b_m)\in\Za[i]^m$, define the \textit{cylinder} of level $m$ based on $\bfb$ as the set
\[
\clC_m(\bfb):=\{z\in \mfF: a_1(z)=b_1,\ldots, a_m(z)=b_m\}
\]
and denote its interior by $\clCc_m(\bfb)$.  Let $\omfF$ be the closure of $\mfF$.  Define the space of sequences
\[
\overline{\sfR} :=
\left\{ \sanu \in \Za[i]^{\Na}: \forall n\in\Na \quad \clCc_n(a_1,\ldots, a_n)\neq \vac \right\}
\]
and let $\overline{\Lambda}:\overline{\sfR}\to \omfF$ be the mapping acting on each $\bfa=\sanu\in\overline{\sfR}$ as follows:
\[
\overline{\Lambda}(\bfa)
:=
\lim_{n\to\infty} [0;a_1,\ldots, a_n]_{\Cx}.
\]
It will be easy to show that $\overline{\Lambda}$ is well defined (see Section \ref{SECTION:PFTHM1}). Loosely speaking, our first result tells us that we may replace $\Omega$ with $\overline{\sfR}$.

\begin{teo01}\label{TEO:01}
$\overline{\Lambda} \left[\,\overline{\sfR}\, \right] = \omfF\setminus \QU(i)$.
\end{teo01}
The inclusion $\overline{\Lambda} \left[\,\overline{\sfR}\, \right] \subseteq \omfF\setminus \QU(i)$ will be a consequence of the approximation properties of Hurwitz continued fractions (Proposition \ref{Propo2.1}.\ref{Propo2.1.v}). We will prove $\omfF\setminus \QU(i) \subseteq \overline{\Lambda} \left[\,\overline{\sfR}\, \right]$, using an iterative procedure (Procedure \ref{Algor}) whose input is a sequence $\bfa\in \Omega\setminus \overline{\sfR}$ and whose output is a sequence $\bfb\in\overline{\sfR}$ such that $\overline{\Lambda}(\bfb)=[0;a_1,a_2,\ldots]_{\Cx}$. In \cite[Corollary 6.4]{BugGeroHus2023}, $\omfF\setminus \QU(i) \subseteq \overline{\Lambda} \left[\,\overline{\sfR}\, \right]$ was obtained by different means but the argument there fails to provide for any $\bfa=(a_n)_{n\geq 1} \in \Omega\setminus\overline{\sfR}$ a sequence $\bfb\in\overline{\sfR}$ for which $\overline{\Lambda}(\bfb)=[0;a_1,a_2,\ldots]_{\Cx}$. 

Let $\sigma\colon \Za[i]^{\Na}\to \Za[i]^{\Na}$ be the left-shift map; that is, $\sigma(a_1,a_2,a_3,\ldots)=(a_2,a_3,\ldots)$ for all $\sanu\in\Za[i]^{\Na}$. In \cite{Yur1995}, Yuri proved that $(\overline{\sfR},\sigma)$ is a sofic shift realized by a directed graph whose vertices are open prototype sets (see definition in Section \ref{SECTION:HCF}) and with edges labeled by Gaussian integers. Yuri's paper concerns the ergodic theory of a class of dynamical systems defined in compact subsets of $\RE^n$, $n\in\Na$. She considered measures equivalent to the $n$-dimensional Lebesgue measure and, thus, certain Lebesgue-null sets and some subsets of the associated shift space are discarded. These omissions force us to make a more detailed study of the relationship between $(\overline{\sfR}, \sigma)$ and the dynamics of Hurwitz continued fraction when dealing with other problems such as those in Section \ref{SECTION:PFTHM1} and Section \ref{SECTION:TRANSCENDENTAL}.

Recall that a real number $x\in [0,1)$ is normal to base $2$ if the asymptotic frequency of every finite word in $\{0,1\}$ of length $k\in\Na$ in the binary expansion of $x$ is exactly $2^{-k}$. In general, given a numeration system, normal numbers are those for which the asymptotic frequency of every finite word of digits is typical. For the unit interval and regular continued fractions, the notion of \textit{typical} is determined by the Gauss measure of cylinders (see \cite{adler1981construction, AirJacKwiMan2020, becher2019absolutely, kraaikamp2000normal}). In the context of Hurwitz continued fractions, we consider the unique Borel probability $T$-ergodic measure $\muh$ which is equivalent to the Lebesgue measure on $\mfF$ (see \cite{EiItoNak2019, Nak1976}) and we refer to normal numbers as \emph{Hurwitz normal numbers} (see Subsection \ref{SUBSEC:HurwitzNormalNumbers} for a precise definition). As a first application of Theorem \ref{TEO:01} and its proof, we determine how complicated is the set of Hurwitz normal numbers as a Borel set.
\begin{teo01}\label{TEO:02}
The set of Hurwitz normal numbers cannot be expressed as a countable intersection of open sets or as a countable union of closed (or even $G_{\delta}$) sets.
\end{teo01}

In fact, we will establish a more robust result by categorizing the set of normal numbers within the framework of the Borel hierarchy. The Borel hierarchy serves as a natural rank, quantifying the number of quantifiers required to construct a Borel set, commencing from open sets. In some sense, this measures the difficulty of determining  whether a complex number is a Hurwitz normal number. A related result was obtained by Dani and Nogueira when they were studying the image of $\Za[i]^2$ under binary quadratic forms \cite[Propositions 8.5, 8.6]{DanNog2014}.

Theorem \ref{TEO:02} is a consequence of the dynamical structure of $\overline{\sfR}$ with respect to the left-shift map. We prove a weak form of specification introduced by Airey, Jackson, Kwietniak, and Mance in~\cite{AirJacKwiMan2020}.
\begin{teo01}\label{TE:RightFeebleSpec}
The subshift $(\overline{\sfR},\sigma)$ has the right feeble-specification property.
\end{teo01}

The subshift $(\overline{\sfR},\sigma)$ is not a full shift, in other words, $\overline{\sfR}$ is not of the form $\scA^{\Na}$ for some $\scA\subseteq \Za[i]$. Nonetheless, Theorem \ref{TE:RightFeebleSpec} shows that the system $(\overline{\sfR}, \sigma)$ still has a considerable level of freedom in terms of specification. Bowen introduced the specification property in his seminal paper \cite{bowen1971periodic}. Bowen's Specification Theorem asserts that all transitive hyperbolic systems possess the ``specification property'', wherein complete independence is established among finite subsets with bounded gaps. The specification and its variations have proven to be potent tools for unlocking access to other properties, such as the uniqueness of measures of maximal entropy \cite{bowen1971periodic,CliTho2021}, universality \cite{quas2016ergodic}, and the intricate structure of generic points \cite{AirJacKwiMan2020,deka2022borel}.

The second application of Theorem 1.1 is related to transcendental numbers. In \cite{AdaBug2007}, Adamczewski and Bugeaud solved Cobham's Conjecture when they proved that, for any $b\in\Na_{\geq 2}$, a real number whose expansion in base $b$ is non-periodic and generated by finite automaton must be transcendental. They extended their work to continued fractions \cite{AdaBug2005} and, afterward, in \cite{Bug2013}, Bugeaud gave a combinatorial condition on non-periodic continued fraction ensuring its transcendence. This condition can be expressed extremely neatly using the repetition exponent introduced by Bugeaud and Kim \cite{BugKim2019}. Given a finite alphabet $\scA$ and a sequence $\bfa=\sanu$ in $\scA$, for each $n\in\Na$ define
\[
r(n,\bfa)
:=
\min\left\{ m\in\Na:  \exists j\in \{1,\ldots, m-n\} \quad (a_j, \ldots, a_{j+n-1} ) = (a_m, \ldots,  a_{m+n-1}) \right\}.
\]
The \textit{repetition exponent} of $\bfa$ is
\[
\rep(\bfa)
:=
\liminf_{n\to\infty} \frac{r(n,\bfa)}{n}.
\]
Bugeaud proved that for any bounded non-periodic sequence $\bfa=\sanu \in \Na^{\Na}$, if $\rep(\bfa)< \infty$, the continued fraction \eqref{Eq:CF} is transcendental for any $a_0\in\Za$ \cite[Theorem 1.3]{Bug2013}. In \cite{Gero2020PP}, the second author obtained a complex version of this result under the addition assumption $\min_n |a_n|\geq \sqrt{8}$. In this paper, we apply the techniques of the proof of Theorem \ref{TEO:01} to construct a family of transcendental numbers whose partial quotients $\sanu$ are bounded and such that $|a_n|=2$ for infinitely many $n\in\Na$. 

 \begin{teo01}\label{TEO:TransComplexNums}
Let $\mathbf{B}:=(B_n)_{n\geq 1}$ be a bounded sequence in $\Za$ which is not periodic and such that 
\[
\min_{n\in\Na} |B_n|\geq 3
\quad\text{ and }\quad
\rep(\bfB)<\infty.
\]
Then, $\zeta=[0;-2,1+iB_1, -2,1+iB_2,-2,1+iB_3,-2,1+iB_4,\ldots]_{\Cx}$ is transcendental.
\end{teo01}

To sum up, these are our applications:
\begin{enumerate}[label= \arabic*.]
\item The set $\overline{\sfR} :=
\left\{ \sanu \in \Za[i]^{\Na}: \forall n\in\Na \quad \clCc_n(a_1,\ldots, a_n)\neq \vac \right\}$ is closed and dynamically chaotic in the sense that it has a weak type of specification. 

\item The set of Hurwitz normal numbers is classified in the third level of the Borel hierarchy.

\item  We construct a family of transcendental complex numbers with bounded partial quotients.
\end{enumerate}

The structure of the paper is as follows. First, in Section \ref{Section:Notation}, we establish notation and conventions. In Section \ref{SECTION:HCF}, we introduce the basic definitions of Hurwitz continued fractions, we discuss the classification of cylinders and prototype sets, and we give some basic results. In Section \ref{SECTION:PFTHM1} we prove Theorem \ref{TEO:01} and state without proof Proposition \ref{LEM:ALGORITHM}, our main technical result. Section \ref{SECTION:PROOF-BOREL-HIERARCHY} starts with some basic notions on Borel Hierarchy and the feeble-specification property. Later, we show Theorem \ref{TE:RightFeebleSpec} and then Theorem \ref{TEO:Nrm:Pi30:Comp} (which implies Theorem \ref{TEO:02}). In Section \ref{SECTION:TRANSCENDENTAL} we show Theorem \ref{TEO:TransComplexNums}. Lastly, Section \ref{SECTION:PFALGORITHM} contains the proof of Proposition \ref{LEM:ALGORITHM}.

\section{Notation and conventions}\label{Section:Notation}
We gather some frequently used notations.
\begin{enumerate}[label= \roman*.]
    \item By natural numbers $\Na$ we mean the set of positive integers and we write $\Na_0:=\Na\cup \{0\}$.
    \item  The real part of a complex number $z$ is denoted by $\real(z)$, its imaginary part by $\imag(z)$, its complex conjugate by $\overline{z}$, and $\Pm(z):=\min\{|\real(z)|, |\imag(z)|\}$.
    \item If $X,Y$ are sets, $A\subseteq X$, $B\subseteq Y$ and $f:X\to Y$ is any function, we define
    \[
    f[A]:=\{f(a):a\in A\}
    \quad\text{and}\quad
    f^{-1}[B]:=\{a\in A: f(a)\in B\}.
    \]
   
    \item $\Za[i]:=\{a+ib\in \Cx: a,b\in\Za \}$, and $\QU(i):=\{a+ib\in \Cx: a,b\in\Za\in \QU\}$.
    
    \item For any $a\in\Za[i]$, the function $\tau_a:\Cx\to\Cx$ is given by $\tau_{a}(z)=z+a$, $z\in\Cx$. By $\iota$ we mean the complex inversion: $\iota(z)=z^{-1}$ for $z\in \Cx\setminus\{0\}$. 
    \item Let $\Rota,\Mir_1:\Cx\to\Cx$ be given by $\Rota(z)=iz$, $\Mir_1(z)=\overline{z}$ and let $\Di_8$ be the group of isometries of $\Cx$ generated by $\Rota$ and $\Mir_1$. Put $\Mir_2=\Rota\circ\Rota\circ\Mir_1$ and denote by $\MIR$ the subgroup of $\Di_8$ generated by $\Mir_1$ and $\Mir_2$.
    \item Given finitely many functions $h_1, \ldots, h_n$ that are either $\iota$, $\tau_a$ for some $a\in\Za[i]$ or some element in $\Di_8$, we write $h_n\cdots h_1:= h_n\circ \cdots \circ h_1$.
    \item For $z\in\Cx$ and $\rho>0$, we write
    \begin{align*}
        \Dx(z;\rho):=\{w\in \Cx: |z-w|<\rho\}, \quad&\quad
        \Dx(z):=\Dx(z;1), \\
        \overline{\Dx}(z;\rho):=\{w\in \Cx: |z-w| \leq \rho\}, \quad&\quad
    \overline{\Dx}(z) := \overline{\Dx}(z;1), \\
    C(z;\rho):=\{w\in\Cx: |z-w|= \rho \}, 
    \quad &\quad
    C(z):=C(z;1).
    \end{align*}

     \item If $X$ is a topological space and $A\subseteq X$, we let $\Cl_X(A)$ (resp. $\inte_X(A)$ ) denote the closure of $A$ (resp. the interior). We omit the sub-index $X$ if the space is clear. 
    \item The term \textit{word} typically refers to a finite set of ordered symbols, while a \textit{sequence} denotes an infinite one. However, we use both terms interchangeably in certain situations where the set may be finite or infinite. 
    \item We denote the empty word by $\epsilon$.
    \item If $\scA$ is a non-empty set, $m\in \Na$, and $\bfa=(a_1,\ldots, a_m)\in\scA^m$, the \textit{length} of $\bfa$ is $|\bfa|:=m$. We further define $|\epsilon|=0$. If $\bfb=\sabu$ is an infinite or finite word in $\scA$, we write 
    \[
    \bfa\bfb:=(a_1,\ldots, a_m,b_1,b_2,\ldots).
    \]
    \item If $\scA$ is a non-empty set and $\bfa=\sanu\in \scA^{\Na}$, a \textit{factor} of $\bfa$ is a subword that appears in any position and a \textit{prefix} is a factor that appears at the beginning of $\bfa$. For each $m\in\Na$ the prefix of length $m$ of $\bfa$ is $\prefi(\bfa;m):=(a_1,\ldots, a_m)$. 
\end{enumerate}

We implicitly use formulas for the inversion of discs and lines in the complex plane throughout the paper. Namely, for $k\in\RE$, we have
\begin{equation}\label{Eq:FormRecInv}
\iota\left[\{z\in \Cx: \real(z)=k\}\right]
=
C\left( \frac{1}{2k}; \frac{1}{2|k|}\right), 
\quad
\iota\left[ \{z\in \Cx: \imag(z)=k\} \right]
=
C\left( \frac{-i}{2k}; \frac{1}{2|k|}\right)
\end{equation}
Also, if $z_0\in\Cx$ and $\rho>0$ are such that $\rho\neq |z_0|$ then
\begin{equation}\label{Eq:FormCircInv}
\iota\left[ C(z_0;\rho)\right]
=
C\left( \frac{\overline{z_0}}{|z|^2-\rho^2}; \frac{\rho}{|\rho^2 - |z_0|^2|}\right).
\end{equation}

\section{Hurwitz continued fractions}\label{SECTION:HCF}
We start this section by defining Hurwitz continued fractions, the oldest complex continued fraction algorithm; other algorithms have been discussed in \cite{ DanNog2014, EiNakadaNatsui2023, LeV1952, AsmusSchmidt1975}. Let $\lfloor \,\cdot\,\rfloor:\RE\to\Za$ be the usual floor function. For $z\in\Cx$, the \textit{nearest Gaussian integer} $[z]$ to $z$ is
\[
[z]
:=
\left\lfloor \real(z) + \frac{1}{2}\right\rfloor
+
i \left\lfloor \imag(z) + \frac{1}{2}\right\rfloor.
\]
We define
\[
 \mfF=\{z\in \Cx:[z]=0\}.
\]

The \textit{Hurwitz map}, $T:\mfF\to \mfF$, is given by
\[
\forall z\in \Cx
\quad
T(z) = 
\begin{cases}
z^{-1}- [z^{-1}], &\text{ if } z\neq 0,\\
0, &\text{ if } z=0.
\end{cases}
\]
For any $z\in \mfF\setminus \{0\}$, we define $a_0(z)\colon=0$, $a_1(z) \colon= [z^{-1}]$ and, if $T^{n-1}(z)\neq 0$, we put $a_n(z)\colon= a_1(T^{n-1}(z))$ for $n\in\Na$. As it is customary,  the exponent in $T$ denotes iteration. For a general $z\in\Cx$, we define $a_0(z)\colon=[z]$ and $a_n(z)\colon=a_n(z-[z])$, $n\in\Na$. The \textit{Hurwitz continued fraction} (\textit{HCF}) of $z$ is
\begin{equation}\label{Eq:DefHCF}
[a_0(z);a_1(z),a_2(z),\ldots]_{\Cx}\colon= 
a_0(z) + \cfrac{1}{a_1(z) + \cfrac{1}{a_2(z) + \cfrac{1}{\ddots}}}.
\end{equation}
It was shown by Hurwitz \cite{Hur1887} that $(a_n(z))_{n\geq 0}$ is infinite if and only if $z\in \Cx\setminus \QU(i)$. In this case, equation \eqref{Eq:DefHCF} can be interpreted as a limit. In the following, we restrict our attention to $z\in\mfF$. For such $z$ 
\begin{equation}
    \label{eq:Hmap}
  [0;a_1(T(z)),a_2(T(z)),\ldots]_{\Cx}=[0;a_2(z),a_3(z)\ldots]_{\Cx}. 
\end{equation}

The terms of the sequence $(a_n(z))_{n\geq 1}$ are the \textit{partial quotients} of $z$. Define the sequences of Gaussian integers $\sepn$, $\seqn$ by:
\begin{equation}\label{Eq:pnqn}
\begin{pmatrix}
p_{-1} & p_{-2} \\
q_{-1} & q_{-2}
\end{pmatrix}
=
\begin{pmatrix}
1 & 0 \\
0 & 1
\end{pmatrix}
\qquad\text{ and }\qquad
\forall n\in\Na_0
\quad
\begin{pmatrix}
p_n\\
q_n
\end{pmatrix}
=
\begin{pmatrix}
p_{n-1} & p_{n-2} \\
q_{n-1} & q_{n-2}
\end{pmatrix}
\begin{pmatrix}
a_n\\
1
\end{pmatrix}.
\end{equation}
The general theory of continued fractions (see, for example, \cite[Section 2]{Khi1997}) tells us that 
\[
\forall n\in\Na
\quad
\frac{p_n}{q_n} = [0;a_1(z),a_2(z),\ldots,a_n(z)]_{\Cx}:=\cfrac{1}{a_1(z)+ \cfrac{1}{a_2(z)+\cfrac{1}{\ddots + \cfrac{1}{a_n(z)}}}}.
\] 
The quotients $\frac{p_n}{q_n}$ are the \textit{convergents} of $z$. From \eqref{Eq:pnqn} we get
\begin{equation}\label{Eq:pnqncoprime}
\forall n\in\Na 
\quad
q_np_{n-1} - q_{n-1}p_n=(-1)^n.
\end{equation}

There are complex numbers $\xi$ for which there are more than one Gaussian integer $a_0$ solving
\begin{equation}\label{Eq:BreakingTies}
|\xi - a_0|
=
\min\left\{ |\xi - a |:a\in \Za[i]\right\}.    
\end{equation}
For such $\xi$, we have defined $[\xi]$ as the solution of \eqref{Eq:BreakingTies} with the largest real and imaginary parts. In many situations, the breaking ties criterion is irrelevant for it applies on a set of Hausdorff dimension $1$, hence Lebesgue null. However, it ensures the uniqueness of the HCF expansion. If we are willing to loosen the expansion's uniqueness slightly, Theorem \ref{TEO:01} and its proof imply that it makes no difference what solution of \eqref{Eq:BreakingTies} we choose.

We will use continued fractions which are not necessarily HCFs. We denote them by
\[
[ b_0;b_1, b_2,\ldots ]
:=
b_0 + \cfrac{1}{b_1 + \cfrac{1}{b_2 + \ddots}}.
\]
We reserve the notation $[b_0;b_1,b_2,\ldots]_{\Cx}$ for HCFs. 

Next, we recall some elementary facts about HCFs. Define $\scD:=\Za[i]\setminus \{0,1,i,-1,-i \}$.
\begin{propo01}\label{Propo2.1}
Take $z=[0;a_1,a_2,\ldots]_{\Cx}\in\mfF\setminus\QU(i)$ and $\sepn$, $\seqn$ as in \eqref{Eq:pnqn}.
\begin{enumerate}[label= (\roman*)., ref = (\roman*)]
\item \label{Propo2.1.i} \cite[Section I]{Hur1887}  For every $n$, we have $a_n\in\scD$.
\item \label{Propo2.1.ii} \cite[Section I]{Hur1887} The sequence $(|q_n|)_{n\geq 0}$ is strictly increasing. 
\item \label{Propo2.1.iii} \cite[Corollary 5.3]{DanNog2014} If $\psi\colon=\sqrt{\frac{1+\sqrt{5}}{2}}$, then $|q_n|\geq \psi^{n-1}$ for all $n\in\Na$.
\item \cite[Theorem 1]{Lak1973} \label{Propo2.1.v} We have
\[
\forall n\in\Na
\quad
\left| z - \frac{p_n}{q_n}\right| < \frac{1}{|q_n|^2}.
\]
\item \cite[Lemma 4.4]{Gero2020PP} \label{Propo2.1.vii} There is a universal constant $\lambda>0$ such that
\[
\forall n\in\Na
\quad
\left| z - \frac{p_n}{q_n}\right| < \frac{\lambda}{|q_nq_{n+1}|}.
\]
\item \cite[Lemma 2.1.(g)]{HeXio2021-01} \label{Propo2.1.viii} There are constants $c_1,c_2>0$ such that for all $n,k\in\Na$ we have
\[
c_1
\leq 
\frac{|q_{n+k}(a_1,\ldots, a_{n+k})|}{|q_n(a_1,\ldots, a_n)||q_k(a_{n+1},\ldots, a_{n+k})|}
\leq 
c_2.
\]
\end{enumerate}
\end{propo01}
\subsection{HCFs Cylinders}\label{SEC:CYL}
\begin{figure}[ht!]
\begin{center}
\includegraphics[scale=0.65,  trim={6.15cm 20cm 1.0cm 4.5cm},clip]{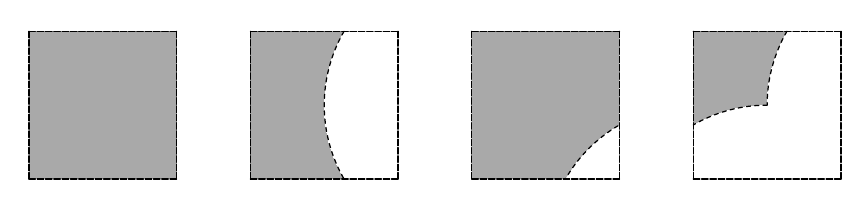}
\caption{ From left to right: $\mfFc$, $\mfFc_1(-2)$, $\mfFc_1(-2+i)$, $\mfFc_1(-1+i)$. \label{Fig-OpenProt}}
\end{center}
\end{figure}

Cylinders and prototype sets are the building blocks of the geometric theory of HCFs. Given $n\in\Na$ and $\bfa=(a_1,\ldots,a_n)\in\scD^n$, the \textit{prototype set} $\mfF_n(\bfa)$ is defined by 
\[
\mfF_n(\bfa)\colon= T^n[\clC_n(\bfa)]. 
\]
We denote the interior (resp. closure) of $\clC_n(\bfa)$ by $\clCc_n(\bfa)$ (resp. $\oclC_n(\bfa)$) and similar conventions apply on $\mfF$. For brevity, by \textit{open cylinders} we mean the interior of a regular cylinder and similarly for prototype sets. Cylinders induce a natural classification of finite words in $\scD$. We say that $\bfa$ is
\begin{multicols}{2}
\begin{enumerate}[label= \roman*.]
\item \textit{valid} if $\clC_n(\bfa)\neq\vac$;
\item \textit{regular} if $\clCc_n(\bfa)\neq \vac$;
\item \textit{irregular} if $\clC_n(\bfa)\neq \vac$ and  $\clCc_n(\bfa)= \vac$;
\item \textit{extremely irregular} if $\#\clC_n(\bfa)=1$.
\end{enumerate}
\end{multicols}
Naturally, $\# A$ represents the number of elements in a set $A$. We denote by $\Omega(n)$ (resp. $\sfR(n)$, $\sfIr(n)$, $\sfEI(n)$) the set of valid (resp. regular, irregular, extremely irregular) words of length $n$. We extend the notions \textit{regular}, \textit{irregular}, and \textit{extremely irregular} to cylinders and prototype sets in an obvious way. It is well known that there are only finitely many prototype sets (see, for example, \cite{EiItoNak2019}). Actually, an open prototype set is one of the following sets or one of their right-angled rotations (see Figure \ref{Fig-OpenProt}):
\[
\mfFc, \quad
\mfFc_1(-2)=\mfFc \setminus \overline{\Dx}(1),
\]
\[
\mfFc_1(-2 + i)=\mfFc \setminus \overline{\Dx}(1-i), \quad
\mfFc_1(-1 + i)=\mfFc \setminus \left(\overline{\Dx}(1)\cup \overline{\Dx}(-i)\right).
\]

A sequence $\bfa=\sanu\in\scD^{\Na}$ is
\begin{enumerate}[label= \roman*.]
\item \textit{valid} if it is the HCF of some $z\in\mfF$;
\item \textit{regular} if it is valid and $(a_1,\ldots, a_n) \in\sfR(n)$ for all $n\in\Na$;
\item \textit{irregular} if it is valid and $(a_1,\ldots, a_n)\in\sfIr(n)$ for some $n\in\Na$;
\item \textit{extremely irregular} if it is valid and $(a_1,\ldots, a_n)\in\sfEI(n)$ for some $n\in\Na$.
\end{enumerate}
We denote by $\Omega$ (resp. $\sfR$, $\sfIr$, $\sfEI$) the set of valid (resp. regular, irregular, extremely irregular) sequences.

A significant part of our calculations is based on an observation stemming from the definitions of cylinders: if $\sanu$ is a sequence in $\scD$, then
\[
\clC_1(a_1) = \mfF \cap \iota\tau_{a_1}[\mfF], 
\quad
\clC_2(a_1,a_2) = \mfF \cap \iota\tau_{a_1}[\mfF]\cap \iota\tau_{a_1}\iota\tau_{a_2}[\mfF];
\]
in general,
\begin{equation}\label{Eq:ObsDefCyls}
\forall n\in\Na
\quad
\clC_n(a_1,\ldots, a_n)
=
\mfF \cap \iota\tau_{a_1}[\mfF] \cap \ldots \cap (\iota\tau_{a_1} \cdots\iota\tau_{a_n})[\mfF].
\end{equation}
Since $\iota$ and $\tau_a$ can be extended to homeomorphisms from the Riemann sphere onto itself, we may replace $\mfF$ and $\clC$ with $\mfFc$ and $\clCc$, respectively, in \eqref{Eq:ObsDefCyls}.  
It follows from \eqref{Eq:ObsDefCyls} that every factor of a regular word (finite or infinite) is again regular. We also deduce the next proposition.

\begin{propo01}\label{PROPO:Recurrencia_mfFn}
Put $\mfF_0(\epsilon):=\mfF$. Every sequence $\sanu$ in $\scD$ satisfies
\[
\forall n\in\Na
\quad
\mfF_n(a_1,\ldots, a_n) 
= 
\tau_{-a_n}\iota\left[\mfF_{n-1}(a_1,\ldots,a_{n-1}) \cap \clC_1(a_n)\right].
\]
\end{propo01}
\subsection{Regular cylinders}
We recall two crucial properties of regular cylinders. The first property tells us that every infinite sequence in $\scD$ with large terms is regular. The second proposition establishes symmetries of the family of regular cylinders.

\begin{propo01}[{\cite[Proposition 4.3]{BugGeroHus2023}}] \label{Propo:Reg01}
If $\bfa=\sanu$ is such that $\min_{n\in \Na} |a_n|\geq \sqrt{8}$, then $\bfa\in\sfR$ and $\mfF_n(a_1,\ldots, a_n)=\mfF$ for all $n\in\Na$.
\end{propo01}

Let $\Rota,\Mir_1:\Cx\to\Cx$ be given by 
\[
\forall z\in \Cx
\quad
\Rota(z):=iz
\quad\text{ and }\quad
\Mir_1(z):=\overline{z}
\]
and call $\Di_8$ the group of isometries they generate. It is clear that $\Di_8$ is the dihedral group of order $8$. Proposition 4.5 and Corollary 4.6 in  \cite{BugGeroHus2023} say that any $f\in \Di_8$ maps open cylinders onto open cylinders of the same level. These results are proven by explicitly determining the image of an open cylinder under $\Rota$ and $\Mir_1$. In particular, if $\Mir_2:=\Rota^2\Mir_1$, $z\mapsto - \overline{z}$, and $\MIR \subseteq \Di_8$ is the subgroup generated by $\Mir_1$ and $\Mir_2$, we have the next result. 
\begin{propo01}\label{Propo:Symm}
If $n\in\Na$, $\bfa=(a_1,\ldots,a_n)\in\sfR(n)$, and $\Mir\in \MIR$, then
\[
\Mir\left[ \clCc_n(\bfa)\right]
=
\clCc_n(\Mir(a_1),\ldots,\Mir(a_n)).
\]
\end{propo01}

\subsection{Irregular cylinders}

\begin{figure}[ht!]
\begin{center}
\includegraphics[scale=0.60,  trim={5.0cm 16.25cm 5.0cm 2.5cm},clip]{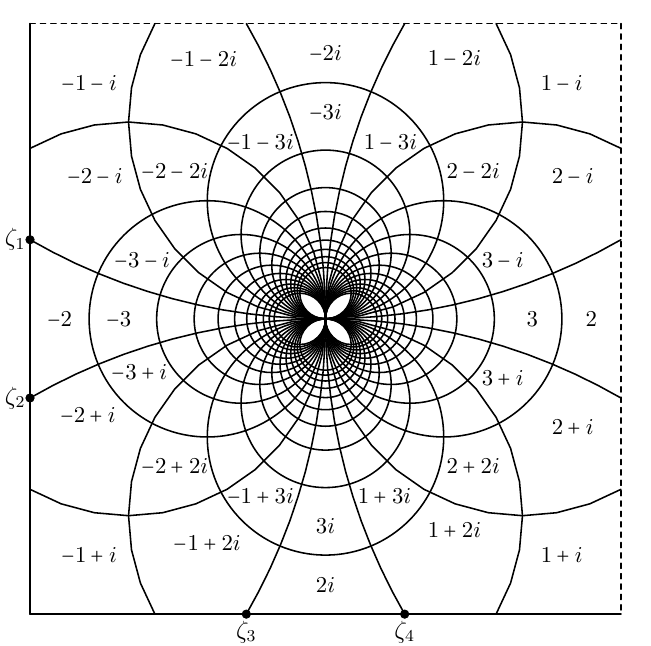}
\caption{The cylinders $\clC_1(a)$, $a\in\scD$, and the numbers $\zeta_1,\zeta_2,\zeta_3,\zeta_4$. \label{Fig-PartmfF}}
\end{center}
\end{figure}

For any pair of complex numbers $z$ and $w$, write
\[
[z,w):=\{z + t(w-z) : t\in [0,1)\}.
\]
Irregular cylinders, and thus prototype sets, exist. Certainly, if
\[
\alpha\colon= \frac{2-\sqrt{3}}{2} \approx 0.13397,
\]
for every $m\in\Za$ with $|m|\geq 2$ we have
\[
\mfF_2(-2,1+mi) = 
\begin{cases}
\left[ -\frac{1}{2} - i\alpha,-\frac{1}{2} + \frac{i}{2}\right), \ \text{\rm if } \  m=2,\\[2ex]
\left[ -\frac{1}{2} - \frac{i}{2},-\frac{1}{2} + i\alpha \right), \ \text{\rm if } \ m=-2, \\[2ex]
\left[ -\frac{1}{2} - \frac{i}{2},-\frac{1}{2} + \frac{i}{2}\right), \  \text{\rm if } \  |m|\geq 3
\end{cases}
\]
(see \cite[Section 4]{BugGeroHus2023}). Additionally, if
\[
\zeta_1 \colon= - \frac{1}{2} + i \alpha, \quad
\zeta_2 \colon= - \frac{1}{2} - i \alpha, \quad
\zeta_3 \colon= - \alpha - \frac{i}{2}, \quad
\zeta_4 \colon=  \alpha - \frac{i}{2},
\]
straightforward calculations yield $a_1(\zeta_1)= -2$, $a_1(\zeta_2)=-2+i$, $a_1(\zeta_3)=2i$, $a_1(\zeta_4)=1 + 2i$ and
\begin{equation}\label{Eq:TsobreZetas}
T(\zeta_1)=\zeta_4, \quad
T(\zeta_2)=\zeta_4, \quad
T(\zeta_3)=\zeta_2, \quad
T(\zeta_4)=\zeta_2.
\end{equation}
As a consequence, we have the HCF expansions
\begin{align*}
\zeta_1 &= [0;   -2, 1+2i, -2+i, 1+2i, -2+i, 1+2i, -2+i, \ldots ]_{\Cx}, \\
\zeta_2 &= [0; -2+i, 1+2i, -2+i, 1+2i, -2+i, 1+2i, -2+i, \ldots ]_{\Cx}, \\
\zeta_3 &= [0;   2i, -2+i, 1+2i, -2+i, 1+2i, -2+i, 1+2i, \ldots ]_{\Cx}, \\
\zeta_4 &= [0; 1+2i, -2+i, 1+2i, -2+i, 1+2i, -2+i, 1+2i, \ldots ]_{\Cx}.
\end{align*}
Using \eqref{Eq:FormCircInv} and \eqref{Eq:ObsDefCyls}, we may show that the sequences $(a_n(\zeta_j))_{n\geq 1}$, $j\in\{1,2,3,4\}$, are extremely irregular. In fact, by carefully tracking the image under $T$ of the boundaries of prototype sets, we may even show that the HCF of a complex number $\xi\in\mfF$ is extremely irregular if and only if $T^n(\xi)=\zeta_4$ for some $n\in\Na$. Figure \ref{Fig-PartmfF} depicts the partition of $\mfF$ induced by the cylinders $\{\clC_1(a):a\in\scD\}$ and the location of $\zeta_1,\zeta_2,\zeta_3,\zeta_4$ in $\mfF$.

\begin{figure}[ht!]
\begin{center}
\includegraphics[scale=0.65,  trim={4.5cm 16.75cm 9.5cm 4.0cm},clip]{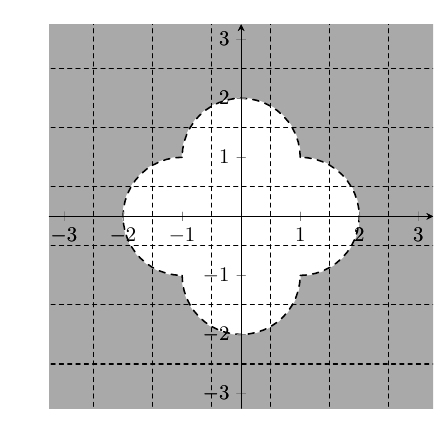}
\caption{ The set $\iota[\mfFc]$. \label{FiginvmfF}}
\end{center}
\end{figure}

\begin{figure}[ht!]
\begin{center}
\includegraphics[scale=0.65,  trim={4cm 17cm 13.5cm 4cm},clip]{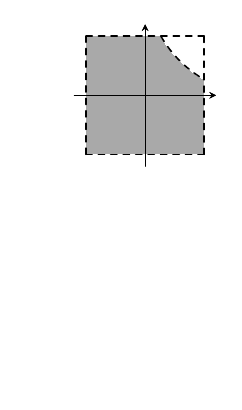}
\includegraphics[scale=0.65,  trim={4.7cm 17cm 10cm 4cm},clip]{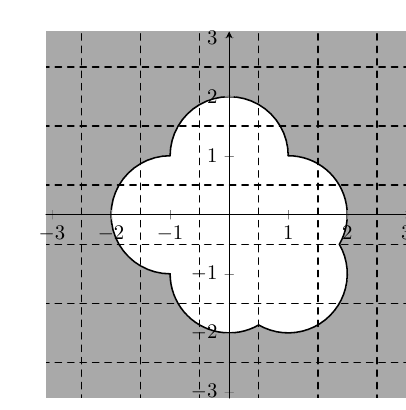}
\caption{ The sets $ \mfFc\setminus \overline{\Dx}(1+i) $ (left) and $\iota[\omfF\setminus \Dx(1+i)]$ (right). \label{FigB5ii}}
\end{center}
\end{figure}

\begin{figure}[ht!]
\begin{center}
\includegraphics[scale=0.65,  trim={4cm 18cm 13.5cm 4cm},clip]{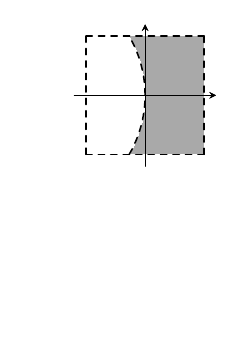}
\includegraphics[scale=0.65,  trim={4.7cm 17cm 10cm 4cm},clip]{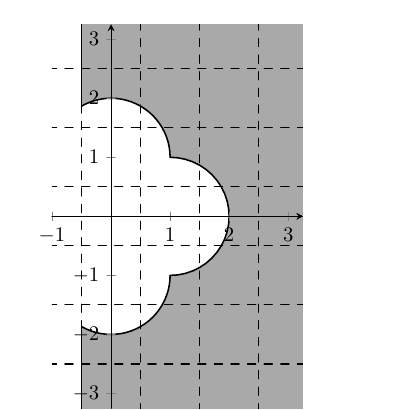}
\caption{ The sets $ \mfFc\setminus \overline{\Dx}(-1)$ (left) and $\iota[ \omfF\setminus \Dx(-1)]$ (right).\label{FigB5iii}}
\end{center}
\end{figure}

\begin{figure}[ht!]
\begin{center}
\includegraphics[scale=0.65,  trim={4cm 18cm 13.5cm 4cm},clip]{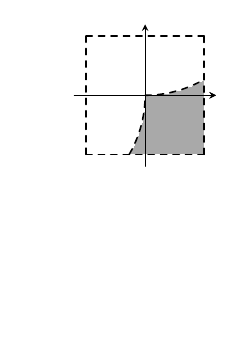}
\includegraphics[scale=0.65,  trim={4.7cm 18cm 10cm 4cm},clip]{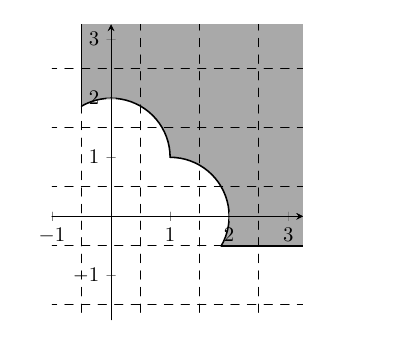}
\caption{ The sets $\mfFc\setminus (\overline{\Dx}(i)\cup \overline{\Dx}(-1)) $ (left) and $\iota[(\omfF\setminus (\Dx(i)\cup \Dx(-1)))]$ (right). \label{FigB5iv}}
\end{center}
\end{figure}

\begin{lem01}\label{Lema:LemaB5}
Let $n\in\Na$ be arbitrary. If $\bfa\in \sfR(n)$ and $b\in\scD$ are such that $\bfa b\in \sfIr(n+1)$, then $\mfFc_n(\bfa)$ is one of the next sets:
\begin{align*}
&\mfFc \setminus \overline{\Dx}(1), \quad
\mfFc \setminus \overline{\Dx}(-i), \quad
\mfFc \setminus \left(\overline{\Dx}(1)\cup \overline{\Dx}(-i)\right),\\
&\mfFc \setminus \left(\overline{\Dx}(-1)\cup \overline{\Dx}(-i)\right), \quad
\mfFc \setminus \left(\overline{\Dx}(1)\cup \overline{\Dx}(i)\right).
\end{align*}
\end{lem01}
\begin{proof}
The proof amounts to checking that every option for non-empty $\mfFc_n(\bfa)$, which are not listed in the statement, satisfies
\[
\bfa b \in\Omega(n+1)
\quad\text{ implies }\quad
\bfa b \in\sfR(n+1).
\]
It follows from the stronger statement
\[
\iota\left[\omfF_n(\bfa)\right] \cap \tau_b\left[ \mfF\right]\neq\vac
\quad\text{ implies }\quad
\iota\left[\mfFc_n(\bfa)\right] \cap \tau_b\left[ \mfFc\right]\neq\vac.
\]
For clarity, we rely on pictures to avoid long yet elementary computations.
\begin{enumerate}[label= \roman*.]
\item If $\mfFc_n(\bfa)=\mfFc$, then $\iota[\mfFc_n(\bfa)]\cap \tau_{b}[\mfFc]\neq \vac$ for all $b\in \scD$ (see Figure \ref{FiginvmfF}).

\item When $\mfFc_n(\bfa) = \mfFc\setminus \overline{\Dx}(1+i)$ we have $\iota[\omfF_n(\bfa)]=\iota[\omfF]\setminus \Dx(1-i)$. We conclude from Figure \ref{FigB5ii} that $\iota\left[\omfF_n(\bfa) \right] \cap \tau_{b}[\mfF] \neq \vac$ if and only if $b\neq 1-i$, in which case $\iota\left[\mfFc_n(\bfa) \right] \cap \tau_{b}[\mfFc ]\neq \vac$. A similar argument holds when $\mfFc_n(\bfa)=\mfFc\setminus i^k\overline{\Dx}(1+i)$ for some $k\in\{2,3,4\}$.
\item If $\mfFc_n(\bfa) = \mfFc \setminus \overline{\Dx}(-1)$, then
\[
\iota\left[\omfF_n(\bfa)\right]
= 
\iota\left[\omfF\right] 
\cap 
\left\{z\in \Cx: \real(z)\geq -\frac{1}{2}\right\}
\]
and, by Figure \ref{FigB5iii}, $\iota\left[\omfF_n(\bfa)\right] \cap \tau_b\left[ \mfF\right] \neq \vac$ if and only if $\real(b)\geq 0$. Figure \ref{FigB5iii} also shows that $\real(b)\geq 0$ yields $\iota[\mfFc_n(\bfa)] \cap \tau_{b}[\mfFc]\neq \vac$. A similar argument holds when $\mfFc_n(\bfa) =  \mfFc \setminus \overline{\Dx}(i)$.

\item When $\mfFc_n(\bfa) = \mfFc \setminus (\overline{\Dx}(i)\cup \overline{\Dx}(-1))$, by Figure \ref{FigB5iv} we have $\iota\left[\omfF_n(\bfa)\right] \cap \tau_b\left[ \mfF\right]\neq\vac$ if and only if $\real(b)\geq 0$ and $\imag(b)\geq 0$. For such $b$ we have $\iota[\mfFc_n(\bfa)] \cap \tau_b[\mfFc]\neq\vac$. 
\end{enumerate}
This shows that if $\bfa\in\sfR(n)$ and $\mfFc_n(\bfa)\neq \vac$ is not one of the sets listed in the lemma, then $\bfa b\in \Omega(n+1)$ implies $\bfa b \in \sfR(n+1)$.
\end{proof}

\section{Proof of Theorem \ref{TEO:01}}\label{SECTION:PFTHM1}
\subsection{The shift space}
\label{section:shiftspace}
Endow $\scD^{\Na}$ with the product topology assuming each factor has the discrete topology. Then $\scD^{\Na}$ can be metricized with the metric $\dist$ given by
\[
\forall \bfa,\bfb\in\scD^{\Na}
\quad
\dist(\bfa,\bfb)
:=
\begin{cases}
0, &\text{ if }\quad  \bfa=\bfb,\\
2^{-\min\{n\in\Na: \prefi(\bfa;n)\neq \prefi(\bfb;n) \}}, &\text{ if }\quad  \bfa\neq \bfb.
\end{cases}
\]
Recall that $\sigma\colon \scD^{\Na} \to \scD^{\Na}$ is the map associating to each $\bfx=\saxu\in \scD^{\Na}$ the sequence $\sigma(\bfx)=(x_j)_{j\geq 2}$.
We link the dynamical systems $(\Omega,\sigma)$ and $(\mfF\setminus \QU(i),T)$ through the map $\Lambda:\Omega\to \mfF\setminus \QU(i)$ given by
\[
\forall \bfa=\sanu\in \Omega
\quad
\Lambda(\bfa) 
:=
[0;a_1,a_2,a_3,\ldots]_{\Cx}.
\]
Let $\bfa\in\Omega$ and $\veps>0$ be arbitrary. 
For $\psi$ as in Proposition \ref{Propo2.1}.\ref{Propo2.1.iii}, take $m\in\Na$ such that $\frac{2}{\psi^{m-1}}<\veps$ and denote by $p_m/q_m$ the $m$th convergent of $\Lambda(\bfa)$. If $\bfb\in \Omega$ is such that $\dist(\bfa,\bfb)< 2^{-m}$, then  $\prefi(\bfa;m)=\prefi(\bfb;m)$ and, by Proposition \ref{Propo2.1}.\ref{Propo2.1.v},
\[
\left| \Lambda(\bfa) - \Lambda(\bfb)\right| 
\leq
\left| \Lambda(\bfa) -  \frac{p_m}{q_m} \right| + \left| \Lambda(\bfb) -  \frac{p_m}{q_m} \right|
\leq \frac{2}{|q_m|^2} 
\leq \frac{2}{ \psi^{m-1} }<\veps.
\]
In other words, $\Lambda$ is uniformly continuous. Call $\Lambda|_{\sfR}$ the restriction of $\Lambda$ to $\sfR$ and let $\overline{\Lambda}$ be the unique continuous extension of $\Lambda|_{\sfR}$ to $\overline{\sfR}:=\Cl_{\scD^{\Na}}(\sfR)$. Note that this definition of $\overline{\sfR}$ is equivalent to the one in the introduction. We now state Proposition \ref{LEM:Algorithm}, which contains the core of Theorem \ref{TEO:01}. We prove it by constructing the sequence $\bfb$ with an iterative process (Procedure \ref{Algor}).
\begin{propo01}\label{LEM:Algorithm}
If $\bfa\in\sfIr$, then there is some $\bfb\in \overline{\sfR}$ such that $\overline{\Lambda}(\bfb)=\Lambda(\bfa)$.
\end{propo01}

The proof of Proposition \ref{LEM:Algorithm} is in Subsection \ref{SUBSECTION:PFPROC}. We now show Theorem \ref{TEO:01}.

\begin{proof}[Proof of Theorem \ref{TEO:01}]
Pick any $\bfa=\sanu\in \overline{\sfR}$ and let $(\bfb^n)_{n\geq 1}$ be a sequence in $\sfR$ converging to $\bfa$. Then,  $\overline{\Lambda}(\bfa)\in \omfF$ because $\overline{\Lambda}$ is continuous and $\overline{\Lambda}(\bfb^n)=\Lambda(\bfb^n) \in \mfF$ for all $n\in\Na$. Next, we show $\overline{\Lambda}(\bfa)\in \omfF\setminus \in\QU(i)$ by obtaining infinitely many coprime pairs $(p,q)\in\Za[i]$ satisfying 
\begin{equation}\label{Eq:Teo1.1:01}
\left| \overline{\Lambda}(\bfa) - \frac{p}{q}\right| 
\leq 
\frac{1}{|q|^2}.
\end{equation}
 Pick $N\in\Na$ and write $p_N=p_N(a_1,\ldots, a_N)$, $q_N=q_N(a_1,\ldots, a_N)$. For all large $n\in\Na$, we have $\prefi(\bfb^n;N)=\prefi(\bfa;N)$ and, by Proposition \ref{Propo2.1}.\ref{Propo2.1.v}, 
\[
\left| \overline{\Lambda}(\bfa) - \frac{p_N}{q_N}\right|
\leq 
\left| \overline{\Lambda}(\bfa) - \overline{\Lambda}(\bfb^n)\right|
+
\left| \overline{\Lambda}(\bfb^n) - \frac{p_N}{q_N}\right| 
\leq
\left| \overline{\Lambda}(\bfa) - \overline{\Lambda}(\bfb^n)\right|
+
\frac{1}{|q_N|^2}.
\]

Letting $n\to\infty$ and using the continuity of $\overline{\Lambda}$, we see that $(p,q)=(p_N, q_N)$ solves \eqref{Eq:Teo1.1:01}. Furthermore, by Proposition \ref{Propo2.1}.\ref{Propo2.1.ii} and \eqref{Eq:pnqncoprime}, the pairs $(p_N,q_N)$, $N\in\Na$, are coprime and different. Therefore, $\overline{\Lambda}(\bfa)\in \omfF\setminus\QU(i)$.

Now consider $z\in \omfF \setminus \QU(i)$. Assume that $z=[0;a_1,a_2,\ldots]\in\mfF$. If $\bfa=\sanu\in \sfR$, then $z=\overline{\Lambda}(\bfa)$. When $\bfa\in\sfIr$, the sequence $\bfb$ from Proposition \ref{LEM:Algorithm} verifies $\Lambda(\bfa)=\overline{\Lambda}(\bfb)$. When $z\not\in\mfF$, we have either $\real(z)= \tfrac{1}{2}$ or $\imag(z)= \tfrac{1}{2}$. If $\imag(z)= \frac{1}{2}$, then $\overline{z}=\Mir_1(z)\in\mfF$ and there is some $\bfb=\sabu\in \overline{\sfR}$ such that $\overline{z}=\overline{\Lambda}(\bfb)$. Thus, by Proposition \ref{Propo:Symm}, we have $\bfc:=(\Mir_1(b_j))_{j\geq 1}\in \overline{\sfR}$ and $\overline{\Lambda}(\bfc)=z$. A similar argument replacing $\Mir_1$ with $\Mir_2$ works when $\real(z)=\frac{1}{2}$.
\end{proof}

\begin{rema}
    Theorem 5.3 in \cite{BugGeroHus2023} implies that $\overline{\Lambda}:\overline{\sfR}\to \omfF$ is at most six-to-one. However, in Section \ref{SECTION:PROOF-BOREL-HIERARCHY} we describe a large set where it is one-to-one.
\end{rema}

\subsection{The construction of $\bfb$}\label{Subsec:Algorithm}

In this subsection, we provide a procedure that gives the sequence $\bfb$ of Proposition \ref{LEM:Algorithm}.

For each $a\in \scD$ such that
\[
a\in\{1+im, -1+im, m+i, m-i\} 
\text{ for some } m\in\Za, \; |m|\geq 2,
\]
define 
\[
S(a)\colon=
\begin{cases}
im, &\text{ if }\quad a\in\{ 1+ im, -1+ im\}, \\
m, &\text{ if }\quad a\in\{ m + i, m - i\}.
\end{cases}
\]
\begin{algor} \label{Algor}
To any $\bfa\in\sfIr$ we associate $\bfb\in\scD^{\Na}$ as in Proposition \ref{LEM:Algorithm}.
\begin{itemize}
\item[•] \verb+Input.+ A sequence $\bfa=\sanu \in \sfIr$.

\item[•] \verb+Output.+ A sequence $\bfb=(b_j)_{j\geq 1} \in \overline{\sfR}$ such that $\overline{\Lambda}(\bfb)=\Lambda(\bfa)$.
\end{itemize}

\begin{enumerate}[label= \textbf{Step }\arabic*.]
\item \label{ProcessStep1} Put $N=0$ and define $\bfb_N\colon= (b_{N}(j))_{j\geq 1}$ by 
\[
\forall j\in\Na
\quad
b_N(j):=a_j.
\]
\item \label{ProcessStep2} If possible, pick $j_N\in \Na$ such that 
\begin{equation}\label{Eq:Algorithm}
(b_N(1),\ldots, b_N(j_N)) \in \sfR(j_N) 
\quad\text{ and }\quad
(b_N(1),\ldots, b_N(j_N+1) ) \not\in\sfR(j_N + 1).
\end{equation}
Define the sequence $\bfb_{N+1}=(b_{N+1}(l))_{l\geq 1}$ as follows:
\begin{enumerate}[label= \roman*.]
\item For $l\in \{1,\ldots,j_N\}$, write 
\[
b_{N+1}(l):=b_{N}(l).
\]
\item For $l=j_N+1$, write
\[
b_{N+1}(j_{N}+1) := S\left( b_{N}(j_{N}+1)\right).
\]
\item For $l\in \{j_{N}+1+n:n\in\Na\}$, let $m\in\Za$ be such that $b_{N}(j_N+1) \in \{1+im, -1+im, m+i, m-i\}$. 
\begin{itemize}
\item[•] If $b_N(j_N+1)\in \{m+i, m-i\}$ for some $m\in\Za$, put
\[
\forall n\in\Na
\quad
b_{N+1}(j_N+1+n):=\Mir_1\left( b_{N}(j_N+1+n)\right).
\]
\item[•] If $b_{N}(j_N+1) \in \{1+im, -1+im\}$, then
\[
\forall n\in\Na
\quad
b_{N+1}(j_N+1+n):=\Mir_2\left( b_{N}(j_N+1+n)\right).
\]
\end{itemize}
\end{enumerate}
If there is no $j_N\in\Na$ such that \eqref{Eq:Algorithm} holds, define $\bfb_{N+1}\colon= \bfb_N$. Update $N:= N+1$.

\item \label{ProcessStep3} Repeat \ref{ProcessStep2}
\item \label{ProcessStep4} Take $\bfb:=\displaystyle\lim_{N\to\infty} \bfb_N$. 
\end{enumerate}
\end{algor}
The reader might wonder whether this process can be performed and, if so, whether the output $\bfb$ has the desired properties. We show in Section \ref{SECTION:PFALGORITHM} the convergence of $(\bfb_N)_{N\geq 1}$, $\bfb\in \overline{\sfR}$, and $\overline{\Lambda}(\bfb)=\Lambda(\bfa)$. 
\subsection{Examples}
For clarity, we illustrate the previous procedure with a couple of examples.
\begin{ex01}
\label{Par:Ex_1}
Let $(B_n)_{n\geq 1}$ be a sequence in $\Za$ such that $|B_n|\geq 3$ for all $n\in\Na$ and let $\bfa=\sanu$ be given by
\[
\bfa = (-2, 1 + iB_1, -2, 1 + iB_2, -2, 1 + iB_3, -2, 1 + iB_4, \ldots).
\]
First, let us show that $\bfa\in\sfIr$. Since $|B_1|\geq 3$, we have
\[
\mfF_2(-2, 1+iB_1) 
=
\left[ -\frac{1}{2} - \frac{i}{2}, -\frac{1}{2} + \frac{i}{2}\right),
\]
so $(-2, 1+iB_1)\in\sfIr(2)$ and $[\zeta_2,\zeta_1]=\oclC_1(-2)\cap \omfF_2(-2, 1+iB_1)\subseteq \mfF_2(-2, 1+iB_1)$ (see Figure \ref{Fig-PartmfF}), which implies 
\[
\oclC_3(-2, 1+iB_1,-2)\subseteq \clC_2(-2, 1+iB_1).
\]
Moreover, from $\mfF_3(-2, 1+iB_1,-2)=\mfF\cap C(1)$, we conclude
\begin{align*}
\iota\left[ \omfF_3 (-2, 1+iB_1, -2)\right] \cap \tau_{1+iB_2}[\omfF]
&=
\left[ 1 + i \left( B_2 - \frac{1}{2}\right), 1 + i \left( B_2 + \frac{1}{2}\right)\right] \\
&\subseteq 
\iota \left[ \mfF_3(-2, 1+iB_1, -2 ) \right].
\end{align*}
and hence, by \eqref{Eq:ObsDefCyls},
\[
\oclC_4(-2, 1+iB_1, -2, 1+iB_2)\subseteq \clC_3(-2, 1+iB_1, -2).
\]
We may continue this way to show 
\[
\forall n\in\Na \quad
\oclC_{n+1}(a_1,\ldots, a_{n+1})  \subseteq \clC_n(a_1,\ldots,a_n).
\]
Since $|\clC_n(a_1,\ldots, a_n)|\to 0$ as $n\to\infty$, we can single out a complex number $z$ satisfying
\[
\{z\} =  \bigcap_{n\in\Na} \clC_n(a_1,\ldots,a_n).
\]
As a consequence, $z=[0;a_1,a_2,\ldots]_{\Cx}$ and $\bfa\in\sfIr$, because $(-2, 1+iB_1)\in\sfIr(2)$. Procedure \ref{Algor} applied to $\bfa$ gives us the sequence 
\[
\bfb= (b_j)_{j\geq 1}= (-2, iB_1, 2, iB_2, -2, iB_3, 2, iB_4, -2, iB_5, 2, iB_6, \ldots),
\]
To see it, we work out the first two iterations.
\begin{enumerate}[label= \roman*.]
\item  Take $N=0$ and $\bfb_0\colon=\bfa$. Then, $j_0=1$ and
\[
b_1(1) := a_1 =-2, \quad
b_1(2) := S(a_2) = i B_1,\quad
\forall l\in \Na_{\geq 3} \;\, b_1(l) := \Mir_2(a_{l}),
\]
which means
\[
\bfb_1
=
(-2, iB_1, 2, -1+iB_2, 2, -1+iB_3, 2, -1+iB_4, \ldots).
\]
Add $1$ to $N$.
\item For $N=1$, we have $j_1=3$ because 
\[
(-2, iB_1, 2)\in\sfR(3)
\quad\text{ and }\quad
(-2, iB_1, 2, -1+iB_2)\not\in\sfR(4).
\]
Note that $(-2, iB_1, 2, -1+iB_2)\not\in\Omega(4)$. We define
\begin{align*}
\forall l\in\{1,2,3\} \quad b_2(l) &:= b_1(l), \\
b_2(4) &:= S(b_1(3)) = i B_2,\\
\forall l\in \Na_{\geq 5} \quad b_2(l)&:= \Mir_2(b_1(l)) = \Mir_2^2(a_l)=a_l,
\end{align*}
which implies
\[
\bfb_2
=
(-2, iB_1, 2, iB_2, -2, 1 + iB_3, -2, 1 + iB_4, \ldots).
\]
\end{enumerate}
We may replace the condition $\min\{|B_n|:n\in\Na\}=3$ by the optimal condition $\min\{|B_n|:n\in\Na\}=2$ at the expense of notational complications.
\end{ex01}

\begin{ex01}
\label{Par:Ex_2} 
Procedure \ref{Algor} applied on the extremely irregular sequence
\[
\bfa= (-2, 1+2i, -2+i, 1+2i, -2+i, 1+2i, -2+i,\ldots),
\]
where the word $(1+2i, -2+i, 1+2i, -2+i)$ is repeated, gives
\[
\bfb=(-2, 2i, 2, -2i,-2, 2i, 2, -2i,\ldots).
\]
Observe that $\zeta_1=\overline{\Lambda}(\bfa)$ and that the sequences
\begin{align*}
    &(-2+i, \, -2i, \,2, \,2i, -2, \, -2i, \,2, \,2i,-2, \, -2i, \,2, \,2i, \ldots), \\
    &(-2+i, \, 1-2i, \,-2, \,2i, \,2, \,-2i, \,-2, \,2i, \,2, \,-2i, \ldots)
\end{align*}
also belong to $\overline{\sfR}$ and are mapped to $\zeta_1$ by $\overline{\Lambda}$ (cfr. \cite[Proposition 5.5]{BugGeroHus2023}). 
\end{ex01}

\section{First application: Borel hierarchy of the set of Hurwitz normal numbers}
\label{SECTION:PROOF-BOREL-HIERARCHY}

This section aims to analyze the complexity of the set of Hurwitz normal numbers with respect to the Borel hierarchy employed in descriptive set theory. This hierarchy provides insight into the intricacy for determining whether a complex number is a Hurwitz normal number. The main step of this proof is to obtain dynamical information of the subshift $\overline{\sfR}$ (Theorem \ref{TE:RightFeebleSpec}). Subsequently, we leverage a result from Airey, Jackson, Kwietniak, and Mance to establish Theorem \ref{TEO:Nrm:Pi30:Comp}  below, which implies Theorem \ref{TEO:02}.

\subsection{Hurwitz normal numbers}\label{SUBSEC:HurwitzNormalNumbers}
In \cite{Nak1976}, Nakada proved the existence and uniqueness of a $T$-ergodic Borel probability measure $\muh$ on $\mfF$ that is equivalent to the Lebesgue measure. By virtue of Birkhoff's Ergodic Theorem \cite[Theorem 1.14]{Walters1982}, for every $n\in\Na$ and every $\bfb \in \sfR(n)$, almost every $z\in \mfF$ with respect to $\muh$ (or, equivalently, to Lebesgue measure) satisfies
\begin{equation}\label{Eq:DefHurwitzNormal}
\lim_{n\to\infty} 
\frac{\#\{j\in\{0,\ldots,n-1\} :(a_{j}(z),\ldots, a_{j+n-1}(z)) = \bfb \}}{n}=\muh ( \clC_n(\bfb)).
\end{equation}
In fact, for almost every $z\in\mfF$ the equality \eqref{Eq:DefHurwitzNormal} holds for every finite regular word $\bfb$. We refer to such numbers $z$ as \textit{Hurwitz normal numbers} and denote them by $\Nrm(\muh)$.

As for regular continued fractions or integer base expansions, normality is closely related to equidistribution. Given a sequence $(z_n)_{n\geq 1}$ in $\mfF$, let $\delta_{z_n}$ be the point mass measure on $z_n$ for each $n\in\Na$. We call $(z_n)_{n\geq 1}$ \textit{uniformly distributed}$\pmod{\muh}$ if the sequence of measures $\frac{1}{n}\sum_{j=1}^n \delta_{z_n}$ converges weakly to $\muh$. 
From the Portmanteau Theorem \cite[Theorem 5.25]{Kal2021}, we can derive the next equality:
\[
\Nrm(\muh)=
\left\{ z\in \mfF: (T^n(z))_{n\geq 1} \text{ is uniformly distributed} \pmod\muh\right\}.
\]
In particular, for every $z\in \Nrm(\muh)$ the set $\{T^n(z):n\in\Na_0\}$ is dense in $\mfF$.

\subsection{Borel hierarchy}\label{Subsec:BH}
The statements in this subsection are classical results in descriptive set theory, their proofs can be found in \cite[Chapter 11]{Kec1995}. Let $(X,\eta)$ be a Polish space; that is, a separable completely metrizable topological space. Denote the first uncountable ordinal by $\omega_1$. We define $\Sigma_1^0(X):=\eta$ and for each ordinal $\xi$ with $1\leq \xi< \omega_1$ put
\begin{align*}
\mathbf{\Pi}_{\xi}^0(X) &\colon= \left\{ X\setminus A\colon A\in \mathbf{\Sigma}_{\xi}^0(X)\right\} \text{, and} \\
 \mathbf{\Sigma}_{\xi}^0(X) &\colon= \left\{ \bigcup_{n\in\Na} A_n: \forall n\in\Na \quad\exists \xi_n<\xi \;\left(A_n\in  \mathbf{\Pi}_{\xi_n}^0(X) \right)\right\}.
\end{align*}
In other words, $\Sigma_1^0(X)$ is the collection of open sets, $\Pi_1^0(X)$ the collection of closed sets, $\Sigma_1^0(X)$ the collection of countable unions of closed sets (an example is $\QU$ when $X=\RE$), and so on. This family of collections has a hierarchical structure; that is, for any two ordinals $\xi,\xi'$ satisfying $1\leq \xi< \xi' < \omega_1$ we have
\[
\mathbf{\Pi}_{\xi}^0(X)\cup \mathbf{\Sigma}_{\xi}^0(X)\subseteq \mathbf{\Pi}_{\xi'}^0(X) 
\quad\text{ and }\quad
\mathbf{\Pi}_{\xi}^0(X)\cup \mathbf{\Sigma}_{\xi}^0(X)\subseteq \mathbf{\Sigma}_{\xi'}^0(X).
\]
Furthermore, every Borel subset belongs to the family, that is, 
\[
\scB(X)
=
\bigcup_{\xi< \omega_1} \mathbf{\Pi}_{\xi}^0(X)
=
\bigcup_{\xi< \omega_1} \mathbf{\Sigma}_{\xi}^0(X), 
\]
where $\scB(X)$ is the Borel $\sigma$-algebra of $(X,\eta)$.
This family of collections is known as the \textit{Borel hierarchy}. 

Note that whenever $X$ is uncountable, all the collections $\mathbf{\Pi}_{\xi}^0(X)$ and $\mathbf{\Sigma}_{\xi}^0(X)$ are different.
The Borel hierarchy is a standard measurement of the complexity of Borel subsets, and it is related to the number of quantifiers needed to describe a subset. For this reason, it is important to know which is the first collection where a Borel set appears. For any ordinal $\xi<\omega_1$, we say that a set $A\in \scB(X)$ is $\mathbf{\Pi}_{\xi}^0(X)$\textbf{-hard} if $A\not\in \mathbf{\Sigma}_{\xi}^0(X)$ . We say that $A\in \scB(X)$ is $\mathbf{\Pi}_{\xi}^0(X)$\textbf{-complete} if it belongs to $\mathbf{\Pi}_{\xi}^0(X)$ and it is $\mathbf{\Pi}_{\xi}^0(X)$-hard. Note that if $A$ is $\mathbf{\Pi}_{\xi}^0(X)$-complete, then $A\in \mathbf{\Pi}_{\xi'}^0(X)\cap \mathbf{\Sigma}_{\xi'}^0(X)$ for every countable ordinal $\xi' >\xi$.

We are ready to state the main result of this section. 
\begin{teo01}\label{TEO:Nrm:Pi30:Comp}
The set $\Nrm(\muh)$ is $\mathbf{\Pi}_3^0(\Cx)$-complete.
\end{teo01}
We will use a variant of the specification to prove this result.

\subsection{Feeble specification}
Let $\scA$ be non-empty and at most countable set. We denote non-empty words $\bfu$ on $\scA$ as $\bfu=(u_1,u_2,\ldots)$. In symbolic dynamics, it is customary to represent words without parentheses, i.e., $\bfu = u_1u_2\cdots$, however, we keep them.

For any two non-empty finite words $\bfv=(v_1,\ldots, v_n)$, $\bfw=(w_1,\ldots, w_n)$ of the same length $n\in\Na$, the \textit{normalized Hamming distance} between $\bfv$ and $\bfw$ is
\[
d_H(\bfv,\bfw) = \frac{\#\{ j\in\{1,\ldots,n\} : v_j\neq w_j\}}{n}.
\] 

Endow $\scA^{\Na}$ with the product topology. 
A \textit{subshift} $X\subset\scA^{\Na}$ is a closed subset such that $\sigma[X]=X$. Given a subshift $X$, we denote by $\scL(X)$ the set of all finite words that appear as factors of points in $X$.
 Given $\bfb\in \scL(\scA^{\Na})\setminus\{\epsilon\}$, we define the \textit{symbolic cylinder based on} $\bfb$ as
\[
[\bfb]_X
\colon=
\left\{ \bfa \in X: \prefi(\bfa;|\bfb|)  = \bfb\right\}.
\]
 The subshift $X$ has the \textit{right feeble specification property} if there is a set $\scG\subseteq \scL(X)$ with the following properties.
\begin{enumerate}[label= \roman*.]
\item If $\bfu,\bfv\in\scG$, then $\bfu\bfv\in\scG$.
\item For every $\veps>0$ there exists $N\in\Na$ such that for any $\bfu\in\scG$ and any $\bfv\in\scL(X)$ satisfying $|\bfv|\geq N$ there are some $\bfs',\bfv'\in \scL(\scA^{\Na})$ satisfying
\[
|\bfv'|=|\bfv|, \quad
0\leq |\bfs|\leq \veps |\bfv|, \quad
d_H(\bfv,\bfv')<\veps, \quad 
\bfu\bfs\bfv'\in\scG.
\]
\end{enumerate}
For any $\bfx=(x_1, x_2, x_3,\ldots)\in\scA^{\Na}$ and $\bfw \in \scL(\scA^{\Na})$ define
\[
e(\bfw,\bfx,N)
\colon=
\#\left\{ j\in \{1,\ldots, N\} : (x_j, x_{j+1}, \ldots, x_{j+n-1}) =\bfw\right\}.
\]
The \textit{quasi-regular set} $Q(X)$ of $X$ is the set of words for which the frequency of every finite word exists:
\[
Q(X)
\:=
\left\{\bfx\in X: \forall \bfw\in \scL(X) \quad  \lim_{N\to\infty} \frac{e(\bfw,\bfx,N)}{N} \text{ exists }\right\}.
\]
For a $\sigma$-invariant probability measure $\nu$, we say that a word $\bfx\in X$ is $\nu$\textit{-generic} if 
\[
\forall \bfw\in \scL(\scA^{\Na})
\quad
\lim_{N\to\infty} \frac{e(\bfw,\bfx,N)}{N} = \nu\left([\bfw]_X\right).
\]
We denote by $G_{\nu}$ the set of $\nu$-generic points. Clearly, $Q(X)\subseteq G_{\nu}$. By Birkhoff's ergodic theorem, if $\nu$ is an $\sigma$-ergodic probability measure, then almost every point is $\nu$-generic.

In \cite{AirJacKwiMan2020}, Airey, Jackson, Kwietniak, and Mance developed a tool to determine the level in the Borel hierarchy of the set of generic sets of subshifts with the feeble specification. 
They applied this result to study normal numbers of regular continued fractions and other numeration systems. 
\begin{teo01}[{\cite[Theorem 6]{AirJacKwiMan2020}}]\label{TEO:TeoAJKM}
Let $X$ be a subshift with the right feeble specification property and at least two $\sigma$-invariant Borel probability measures. Let $\nu$ be a shift-invariant Borel probability measure. If $B$ is a Borel set such that $G_{\nu}\subseteq B\subseteq Q(X)$, then $B$ is $\mathbf{\Pi}_3^0(X)$-hard. In particular, $G_{\nu}$ and $Q(X)$ are $\mathbf{\Pi}_3^0(X)$-complete.
\end{teo01}

\subsection{Proof of Theorem \ref{TE:RightFeebleSpec}}

The next lemma follows from \cite[Proposition 8.4]{DanNog2014}. 

\begin{lem01}\label{Prop:PropConcat}
Let $n\in \Na$ and $\bfa=(a_1,\ldots, a_n) \in \sfR(n)$ be arbitrary.
\begin{enumerate}[label= \roman*.]
\item\label{Prop:PropConcat:01} If $\Pm(a_n)\geq 3$ and $\bfb\in \sfR(m)$ for some $m\in\Na$, then $\bfa\bfb\in\sfR(m+n)$.
\item\label{Prop:PropConcat:02} There exists some $b\in\scD$ such that $\Pm(b)\geq 3$ and $\mfF_{n+1}(\bfa b)=\mfF$.
\end{enumerate}
\end{lem01}



\begin{proof}[Proof of Theorem \ref{TE:RightFeebleSpec}]
By definition, $\overline{\sfR}$ is closed in $\scD^{\Na}$. Since factors (subwords) of regular words are also regular, $\sigma[\overline{\sfR}] \subseteq \overline{\sfR}$. If $\bfa\in\overline{\sfR}$ and $b=3+3i$, then $b\bfa\in \overline{\sfR}$ and $\sigma(b\bfa)=\bfa$, so $\overline{\sfR}\subseteq \sigma[\overline{\sfR}]$. Define the set
\[
\scG
\colon= 
\left\{ (c_1,\ldots, c_m) \in \bigcup_{n\in\Na} \sfR(n): \Pm(c_m)\geq 3\right\}.
\]
Let us verify the conditions defining the right feeble specification property.
\begin{enumerate}[label= \roman*.]
\item If $\bfu,\bfv\in\scG$, then $\bfu\bfv\in\scG$ by Lemma \ref{Prop:PropConcat}.\ref{Prop:PropConcat:01}.
 
\item Let $\veps>0$ and $N\in\Na$ such that $1/N<\veps$. Take $\bfv=(v_1,\ldots, v_m)\in\scL(\overline{\sfR})$ with $m\geq N$. Then $\bfv\in \sfR(m)$, so $(v_1,\ldots,v_{m-1})\in \sfR (m-1)$ and, by Lemma \ref{Prop:PropConcat}.\ref{Prop:PropConcat:02}, there exists $a\in\scD$ satisfying
\[
\Pm(a)\geq 3
\quad\text{ and }\quad
\bfv'\colon= (v_1,\ldots,v_{m-1},a)\in \scG.
\]
Letting $\bfs$ be the empty word, for every $\bfu\in\scG$ we have
\[
|\bfs|\leq \veps |\bfv|, \quad
|\bfv'|=|\bfv|,\quad 
d_H(\bfv,\bfv')<\veps \quad
\bfu\bfs\bfv'\in\scG. \; \qedhere
\]
\end{enumerate}
\end{proof}
\subsection{Proof of Theorem \ref{TEO:Nrm:Pi30:Comp}}
We split the proof of Theorem \ref{TEO:Nrm:Pi30:Comp} into Lemma \ref{LEM:Nrm:Pi30:01} (lower bound) and Lemma \ref{LEM:Nrm:Pi30:02} (upper bound).

 Let $X,Y$ be Polish spaces. We call $A\subseteq X$ \textit{Wadge reducible} to $B\subseteq Y$ if there exists a continuous function $f:X\to Y$ (called a \textit{Wadge reduction}) such that $f^{-1}[B]=A$. For more on Wadge reductions, see \cite[Section 21.E]{Kec1995}.
 
\begin{rema}\label{Rema:WadgeReduction}
   Using basic properties of continuous functions, one can show that if a set $A$ is $\mathbf{\Pi}_{3}^0(X)$-hard and $A\subset X$ is Wadge reducible to $B\subset Y$, then $B$ is $\mathbf{\Pi}_{3}^0(Y)$-hard.  
\end{rema}

\begin{lem01}\label{LEM:Nrm:Pi30:01}
The set $\Nrm(\muh)$ is $\mathbf{\Pi}_3^0(\Cx)$-hard.
\end{lem01}
\begin{proof}

 By Proposition \ref{Propo:Reg01}, any constant sequence $\bfa=(a, a,\ldots)\in\scD^{\Na}$ with $|a|\geq \sqrt{8}$ belongs to $\overline{\sfR}$ and the point mass measure based on $\bfa$ is $\sigma$-invariant. This implies that $(\overline{\sfR},\sigma)$ has infinitely many $\sigma$-invariant measures.

Let us construct the $\sigma$-invariant Borel probability measure $\nu$ on $\overline{\sfR}$ on which we apply Theorem \ref{TEO:TeoAJKM}. We define the set
\[
\scI
\colon= 
\bigcap_{n\in\Na} \bigcup_{\bfa\in \sfR(n)} \clCc_n(\bfa)\subseteq \mfF.
\]
It is clear that $\scI\subset \Cx$ is a Borel subset.
We will use the continuous functions $\Lambda$ and $\overline{\Lambda}\colon \overline{\sfR}\to \Cx$ defined in Subsection \ref{section:shiftspace}.


We define the Borel probability measure $\nu$ on $\overline{\sfR}$ as follows.
\begin{enumerate}[label= \roman*.]
\item We set $\nu(B):=\muh(\overline{\Lambda}[B])$ for every Borel subset $B\subset \overline{\Lambda}\,^{-1}[\scI].$

\item $\nu\left( \overline{\sfR} \setminus \overline{\Lambda}\,^{-1}[\scI] \right)=0$.
\end{enumerate}

Since $\muh$ is $T$-invariant and $\overline{\Lambda}\circ \sigma = T\circ \overline{\Lambda}$ (Equation \eqref{eq:Hmap}), we have that $\nu$ is $\sigma$-invariant. Theorem \ref{TEO:TeoAJKM} implies that $G_{\nu}\subseteq \overline{\sfR}$ is $\Pi_3^0(\overline{\sfR})$-complete.


To finish the proof, we must translate this conclusion to the complex plane using $\overline{\Lambda}$ as a Wadge reduction. First, note that if $z$ belongs to the boundary of some prototype set, then $T(z)$ also belongs to the boundary of some prototype set (the proof of this assertion amounts to applying the complex inversion to a finite collection of line segments and arches). Since for each $z\in\Nrm(\muh)$ the set $\{T^n(z):n\in\Na_{\geq 0}\}$ is dense in $\mfF$, we have $\Nrm(\muh)\subseteq \scI$. Hence, the definitions of $G_{\nu}$ and $\Nrm(\muh)$ yield
\[
G_{\nu}=\overline{\Lambda}\,^{-1}[\Nrm(\muh)].
\]
Therefore $G_{\nu}\subseteq \overline{\sfR}$ is Wadge reducible to $\Nrm(\muh)\subseteq \Cx$. Using Remark \ref{Rema:WadgeReduction} we conclude that $\Nrm(\muh)$ is $\mathbf{\Pi}_3^0(\Cx)$-hard. 
\end{proof}

\begin{lem01}\label{LEM:Nrm:Pi30:02}
The set $\Nrm(\muh)$ is $\mathbf{\Pi}_3^0(\Cx)$.
\end{lem01}

\begin{proof}
Let $(\veps_k)_{k\geq 1}$ be a strictly decreasing sequence of positive numbers converging to $0$. Given $A\subseteq \mfF$, let $\chi_A$ be its indicator function. For $M,n,k\in\Na$ and $\bfb\in \sfR(n)$, define the sets
\begin{align*}
E_{n}(\bfb,M,k) &:= \left\{z = [0;a_1,a_2,\ldots]_{\Cx} : \left|\frac{1}{M} \sum_{j=0}^{M-1} \chi_{\clC_n(\bfb)}\left( T^j(z)\right) - \muh(\clC_n(\bfa))\right| \leq \veps_k\right\}, \\
F_{n} (\bfb,M,k) &:= \Cl_{\Cx}(E_{n}(\bfb,M,k))\text{, and}\\
G_{n} (\bfb,M,k) &:= \inte_{\Cx}(E_{n}(\bfb,M,k)).
\end{align*}
Note that each $E_{n}(\bfb,M,k)$ is a finite union of cylinders of level $n+M-1$. The definition of $\Nrm(\muh)$ gives
\[
\Nrm(\muh) 
=
\bigcap_{n\in\Na} \bigcap_{\bfb\in\sfR(n)} \bigcap_{k\in\Na} \bigcup_{N\in \Na} \bigcap_{M\in\Na_{\geq N}} E_{n}(\bfb,M,k). \nonumber\\
\]
Put
\begin{align*}
\mathfrak{L} &\colon= \bigcap_{n\in\Na} \bigcap_{\bfb \in\sfR(n)} \bigcap_{k\in\Na} \bigcup_{N\in \Na} \bigcap_{M\in\Na_{\geq N}} G_{n}(\bfb,M,k), \nonumber\\
\mathfrak{M} &\colon= \bigcap_{n\in\Na} \bigcap_{\bfb \in\sfR(n)} \bigcap_{k\in\Na} \bigcup_{N\in \Na} \bigcap_{M\in\Na_{\geq N}} F_{n}(\bfb,M,k), \nonumber
\end{align*}
so $\mathfrak{L} \subseteq \Nrm(\muh) \subseteq \mathfrak{M}$. Note that $\bigcup_{N\in \Na} \bigcap_{M\in\Na_{\geq N}} F_{n}(\bfb,M,k)\in \mathbf{\Sigma}_2^0(\Cx)$, this implies that $\mathfrak{M}\in \mathbf{\Pi}_3^0(\Cx)$.  
We will show that $\mathfrak{M}\subseteq \mathfrak{L}$ to conclude the lemma. Let $z\in \mathfrak{M}$ be arbitrary and assume that 
\begin{equation}\label{Eq:UpperBound:01}
    z\in \bigcap_{n\in \Na} \clCc_n\left(a_1(z),\ldots, a_n(z)\right).
\end{equation}
Let $n,k\in \Na$, $\bfb\in \sfR(n)$ and consider $N\in\Na$ with
\[
\forall M\in \Na_{\geq N} 
\quad
z\in F_n(\bfb, M,k). 
\]
Take any $M\in\Na_{\geq N}$ and let $\bfd\in \sfR(M+n)$ be such that 
\[
\clC_{M+n}(\bfd)\subseteq E_n(\bfb,M,k)
\quad\text{ and }\quad
z\in \oclC_{M+n}(\bfd).
\]
Since different regular cylinders of the same level have disjoint interiors, from \eqref{Eq:UpperBound:01} we get $$\bfd = \left(a_1(z),\ldots, a_{M+n}(z)\right),$$ so $z\in G_n(\bfb,M,k)$. We conclude $\mathfrak{M} \subseteq \mathfrak{L}$, hence $\Nrm(\muh)=\mathfrak{M}$, if \eqref{Eq:UpperBound:01} holds for all $z\in\mathfrak{M}$.

Lastly, we show that \eqref{Eq:UpperBound:01} is indeed true for all $z\in\mathfrak{M}$. To this end, pick any $z\in\mathfrak{M}$. It suffices to prove that $\{T^n(z):n\in\Na_{\geq 0}\}$ is dense in $\mfF$, which in turns follows from
\begin{equation}\label{Eq:UpperBound:02}
\forall n\in\Na
\quad
\forall \bfb\in \sfR(n)
\quad
\{T^n(z):n\in\Na_{\geq 0}\}
\cap \clCc_n(\bfb)\neq \vac.
\end{equation}
Take any $n\in\Na$ and $\bfb\in\sfR(n)$. Let $d\in\scD$ be such that $\Pm(d)\geq 3$ and $c\in\scD$ such that
\[
\bfb c \in \sfR(n+1), 
\quad\text{ and }\quad
\oclC_{n+1}( \bfb c)  \subseteq \clCc_n(\bfb)
\]
(see Lemma \ref{LEM:ALGORITHM} below). Since $d\bfb c\in\sfR(n+2)$ by Lemma \ref{Prop:PropConcat}. \ref{Prop:PropConcat:01}, we may choose $k\in\Na$ so large that
\[
\veps_k < \frac{\muh(\clC_{n+2}(d\bfb c))}{4}
\]
and $N\in\Na$ with
\[
\forall M\in\Na_{\geq N} \quad z\in F_{n+2}(d \bfb c, M, k).
\]
Take $\bfd\in \sfR(M+n-1)$ satisfying
\[
\clC_{M+n-1}(\bfd)\subseteq E_{n+2}(d \bfb c,M,k)
\quad\text{ and }\quad
z\in\oclC_{M+n-1}(\bfd).
\]
The choice of $k$ and the definition of $E_{n+1}(\bfb c, M,k)$ guaranty the existence of some $r\in\{n+2,\ldots, M+n-1\}$ satisfying
\[
(d_{r-n-1},\ldots, d_{r-1},d_r)
=
(d, b_1,\ldots, b_n,c).
\]
Let us show that
\begin{equation}\label{Eq:DoubleContention}
\oclC_{M+n-1}(\bfd)
\subseteq
\oclC_{r}(d_1, \ldots, d_r) 
\subseteq 
\clCc_{r-1}(d_1, \ldots, d_{r-1}).
\end{equation}
The first contention is trivial. From 
\begin{align*}
\clC_{r}(\bfd) &= \iota \tau_{d_1}\cdots \iota \tau_{d_r}[\mfF_{r}(d_1,\ldots,d_r)],  \\  
\clC_{r-1}(d_1,\ldots, d_{r-1} ) &= \iota \tau_{d_1}\cdots \iota \tau_{d_{r-1} }[\mfF_{r-1}(d_1,\ldots,d_{r-1})].
\end{align*}
and the fact that $\iota$ and the functions $\tau_a$, $a\in\scD$, are homeomorphisms from the Riemann sphere onto itself, the second contention in \eqref{Eq:DoubleContention} is equivalent to 
\begin{equation}\label{Eq:DoubleContention:02}
    \iota \tau_c\left[ \omfF_{r}(\bfd' d \bfb c) \right]
    \subseteq 
    \mfFc_{r}(\bfd' d \bfb)
    \;\text{ for }\; \bfd'=(d_1,\ldots, d_{r-n-2}).
\end{equation}
Since $\Pm(d)\geq 3$, we may use Proposition \ref{PROPO:Recurrencia_mfFn} to show that 
\[
\mfF_r(\bfd'd\bfb c)=\mfF_r(\bfb c) 
\; \text{ and }\;
\mfF_{r-1}(\bfd' d \bfb)=\mfF_n(\bfb). 
\]
Therefore, the second contention is \eqref{Eq:DoubleContention} is equivalent to $\iota\tau_c[\omfF_{n+1}(\bfb c)]\subseteq \mfFc_n(\bfb)$, which is equivalent to $\oclC_{n+1}(\bfb c)\subseteq \clCc_n(\bfb)$. We conclude $z\in \clCc_n(d_1,\ldots, d_{r-1})$ so $T^{r-n-1}(z)\in \clCc_n(\bfb)$ and \eqref{Eq:UpperBound:02} follows.
\end{proof}
\section{Second application: construction of transcendental numbers}\label{SECTION:TRANSCENDENTAL}
Let $\scA$ be a finite alphabet and $\bfa=\sanu$ a sequence on $\scA$. Recall that for each $n\in\Na$ we have defined
\[
r(n,\bfa)
:=
\min\left\{ m\in\Na:  \exists j\in \{1,\ldots, m-n\} \quad (a_j, \ldots, a_{j+n-1} ) = (a_m, \ldots,  a_{m+n-1}) \right\}.
\]
The \textit{repetition exponent} of $\bfa$ is
\[
\rep(\bfa)
:=
\liminf_{n\to\infty} \frac{r(n,\bfa)}{n}.
\]

Bugeaud and Kim \cite{BugKim2019} showed that the finiteness of $\rep(\bfa)$ for a given sequence $\bfa$ is equivalent to a combinatorial condition used in \cite{AdaBug2005, Bug2013} to construct real transcendental numbers. 

\begin{lem01}[{\cite[Section 10]{BugKim2019}}]\label{LEM:WUVU}
Let $\scA$ be a finite alphabet and $\bfa= \sanu\in \scA^{\Na}$. If $\rep(\bfa)<\infty$, then there are three sequences $(W_n)_{n\geq 1}$, $(U_n)_{n\geq 1}$, $(V_n)_{n\geq 1}$ of finite words in $\scA$ such that
\begin{enumerate}[label= \arabic*.]
\item For all $n\in \Na$, the word $W_nU_nV_nU_n$ is a prefix of $\bfa$;
\item The sequence $\left( (|W_n| + |V_n|)/|U_n|\right)_{n\geq 1}$ is bounded above;
\item The sequence $(|U_n|)_{n\geq 1}$ is increasing.
\end{enumerate}
\end{lem01}

Theorem \ref{TEO:GeroTransHCF:PP} below is a complex version of Bugeaud's transcendence criterion for non-periodic continued fractions with finite repetition exponent \cite{Bug2013}.
\begin{teo01}[{\cite[Theorem 5.1]{Gero2020PP}}] \label{TEO:GeroTransHCF:PP}
Let $\bfa\in\Omega$ be a bounded and non-periodic sequence such that $\rep (\bfa) < \infty$. Let $(W_n)_{n\geq 1}$, $(U_n)_{n\geq 1}$, and $(V_n)_{n\geq 1}$ be sequences as in Lemma \ref{LEM:WUVU}.
\begin{enumerate}[label= \arabic*.]
\item \label{TEO:GeroTransHCF:PP:01} If $\displaystyle\liminf_{n\to \infty} |W_n|<\infty$, then $\Lambda(\bfa)$ is transcendental.
\item \label{TEO:GeroTransHCF:PP:02} If $\displaystyle\liminf_{n\to \infty} |W_n|=\infty$ and $|a_n|\geq \sqrt{8}$ for all $n\in\Na$, then $\Lambda(\bfa)$ is transcendental.
\end{enumerate}
\end{teo01}

We conjecture that the condition $|a_n|\geq \sqrt{8}$ for all $n\in\Na$ can be removed in Theorem \ref{TEO:GeroTransHCF:PP}.\ref{TEO:GeroTransHCF:PP:02}. The second application of our work, Theorem \ref{TEO:TransComplexNums}, is constructing a family of transcendental numbers overlooked by Theorem \ref{TEO:GeroTransHCF:PP}. This provides evidence supporting our conjecture. A key ingredient of the proof is Schmidt's extension of his celebrated Subspace Theorem to number fields, we only state it for $\QU(i)$. 
For $m\in\Na$ and $\bfz=(z_1,\ldots,z_m)\in\Cx^m$, we write $\|\bfz\|=\max\{ |z_1|,\ldots, |z_m|\}$ and $\overline{\bfz}=(\overline{z_1}, \ldots, \overline{z_m})$.


\begin{teo01}[{\cite[Theorem 3]{Sch1975Subspace}}] \label{TEO:SST}
Let $m\in \Na$. Assume that we have two collections of linearly independent linear forms $\scL_{1}, \ldots, \scL_m$ and $\scM_1, \ldots, \scM_m$ in $m$ variables with real or complex algebraic coefficients. Let $\veps>0$. Then, there are finitely many proper subspaces $T_1, \ldots, T_k$ of $\Cx^m$ containing all the solutions $\bfz=(z_1,\ldots, z_m)\in\Za[i]^m$ of
\[
\left|\scL_1(\bfz) \cdots \scL_m(\bfz)\right|
\left|\scM_1(\overline{\bfz}) \cdots \scM_m(\overline{\bfz})\right|
<  \| \bfz\|^{-\veps}.
\]
\end{teo01}
Recall that if $A_0, A_1,\ldots, A_n$ are variables and $q_0,q_1,\ldots, q_n$ are the denominators of the continued fraction $[ A_0; A_1,A_2,A_3,\ldots, A_n]$, then
\begin{equation}\label{Eq:MirrorFormula}
\frac{q_{n-1}}{q_n}
=
[ 0; A_n, A_{n-1}, \ldots, A_2, A_1 ]
\end{equation}
(see, for example, \cite[Theorem 6]{Khi1997}). Given $\bfB\in\Omega$ as in Theorem \ref{TEO:TransComplexNums}, let $(W_n)_{n\geq 1}$, $(U_n)_{n\geq 1}$, $(V_n)_{n\geq 1}$ be as in Lemma \ref{LEM:WUVU} and define the sequences of non-negative integers $(w_n)_{n\geq 1}$, $(u_n)_{n\geq 1}$, and $(v_n)_{n\geq 1}$ by
\[
\forall n\in\Na
\quad
w_n:=|W_n|, \;
u_n:=|U_n|, \;
v_n:=|V_n|, \;
t_n:=w_n + u_n + v_n.
\]
The following lemma is essentially \cite[Lemma 5.2]{Gero2020PP}.
\begin{lem01}\label{LEM:TRANS_EPSILON}
Let $\zeta, \bfB$ be as in Theorem \ref{TEO:TransComplexNums} and let $(q_n)_{n\geq 1}$ correspond to $\zeta$ as in Equation \eqref{Eq:pnqn}. Then, there exists an $\veps>0$ such that for all $n\in\Na$ we have $\psi^{2u_n} \geq  |q_{2w_n}{q_{2t_n}}|^{\veps}$. 
\end{lem01}

For any two finite words of the same length $A = (a_1, \ldots , a_n)$ and $B = (b_1, \ldots, b_n)$ the word $s(A,B):=(c_1,\ldots,c_{2n})$ is given by $c_{2j-1} = a_j$ and $c_{2j}=b_j$ for all $j\in\{1,\ldots, n\}$; that is,
\[
s(A,B)= (a_1,\, b_1,\, a_2,\, b_2, \ldots, \, a_{n},\, b_n).
\]
Now we show Theorem \ref{TEO:TransComplexNums}. Similar arguments have appeared in \cite{AdaBug2005, Bug2013} for real numbers and in \cite[Section 5.2]{Gero2020PP} for complex numbers. 
\begin{proof}[Proof of Theorem \ref{TEO:TransComplexNums}]
Assume that $\zeta$ is an algebraic number, so $[\QU(i,\zeta): \QU(i)]\geq 3$ (its HCF is not periodic; see, for example, \cite[Theorem 4.2]{DanNog2014}). Let $(W_n)_{n\geq 1}$, $(U_n)_{n\geq 1}$, and $(V_n)_{n\geq 1}$ be the sequences associated to $\bfB$ as in Lemma \ref{LEM:WUVU} and define $(w_n)_{n\geq 1}$, $(u_n)_{n\geq 1}$, and $(v_n)_{n\geq 1}$ as above. We may assume that every $u_n$ is even. Indeed, suppose that $u_n$ is odd for every large $n$. For such $n$, call $x$ be the last letter of $U_n$ and $\widehat{U}_n:=\prefi(U_n, u_n-1)$, so
\[
W_nU_nV_nU_n
=
W_n \, \widehat{U}_n x \, V_n  \widehat{U}_n x
\]
and we can replace $U_n$ with $\widehat{U}_n$ and $V_n$ with $x \,V_n$.

When $\displaystyle\liminf_{n\to\infty} w_n<\infty$, the sequence $(-2, 1+iB_1, -2,1+iB_2, \ldots)$ also has finite repetition exponent and its corresponding sequence $(W_n)_{n\geq 1}$ has a constant subsequence. The transcendence of $\zeta$ then follows from Theorem \ref{TEO:GeroTransHCF:PP}. Suppose that $w_n\to\infty$ when $n\to\infty$. As noted in \cite[p. 1013]{Bug2013}, we do not lose any generality by assuming $B_{w_n}\neq B_{t_n}$ for all $n\in\Na$. By Example \ref{Par:Ex_1}, we have $\bfa\in \sfIr$ and the output $\bfb$ of Procedure \ref{Algor} is
\[
\bfb
=
(b_1,b_2,\ldots)
=
(-2, iB_1, 2, iB_2, -2, iB_3, 2, iB_4, \ldots).
\]
Note that
\begin{equation}\label{Eq:Trans:00}
    \forall n\in\Na
    \quad
    \{(b_1,\ldots, b_n), \; (b_n, \ldots, b_1)\}\subseteq  \sfR(n).    
\end{equation}
For each $n\in\Na$, let $A_n$, $B_n$, $C_n$ be the alternating words 
\[
A_n := \left( -2, \, 2,  \,\ldots, \, (-1)^{w_n}2\right),\quad
B_n := \left( (-1)^{w_n+1}2, \, (-1)^{w_n+2}2, \ldots, \, (-1)^{w_n+u_n}2 \right),
\]
\[
C_n := \left( (-1)^{w_n+1}2, \, (-1)^{w_n+2}2, \ldots, \, (-1)^{w_n+v_n}2\right),
\]

and define
\[
W_n':= s(A_n,W_n), \quad
U_n':= s(B_n,U_n), \quad
V_n':= s(C_n,V_n).
\]
By Lemma \ref{LEM:ALGORITHM} below and $\min_j |B_j|\geq 3$, we have $\mfF_{2k}(\prefi(\bfb,2k))=\mfF$ for all $k\in\Na$. Then, since $u_n$ is even, the periodic sequence $Z_n:=W_n'\, U_n'\,V_n'\,U_n'\,V_n'\,U_n'\cdots$ belongs to $\overline{\sfR}$ and $\zeta_n:=\overline{\Lambda}(Z_n)$ is a quadratic irrational over $\QU(i)$. Let $(p_k)_{k\geq 0}$, $(q_k)_{k\geq 0}$ correspond to $\bfb$ as in \eqref{Eq:pnqn}. One can use the standard proof of Lagrange's theorem on periodic continued fractions (see, for example, \cite[Theorem 4.2]{DanNog2014}) to check that $\zeta_n$ is a solution of the polynomial 
\[
P_n(X) 
:= 
\begin{vmatrix} q_{2w_n-1} & q_{2t_n-1} \\ q_{2w_n} & q_{2t_n} \end{vmatrix} X^2 - \left( \begin{vmatrix} q_{2w_n-1} & p_{2t_n-1} \\ q_{2w_n} & p_{2t_n} \end{vmatrix}  + \begin{vmatrix} p_{2w_n-1} & q_{2t_n-1} \\ p_{2w_n} & q_{2t_n} \end{vmatrix} \right)X + \begin{vmatrix} p_{2w_n-1} & p_{2t_n-1} \\ p_{2w_n} & p_{2t_n}  \end{vmatrix}.
\]
We will use $P_n$ to obtain three collections of linear forms and then we will apply Theorem \ref{TEO:SST} to get a contradiction. From $\zeta,\zeta_n\in \oclC_{2tn+2u_n}(b_1, \ldots, b_{2tn+2u_n})$ and Proposition \ref{Propo2.1}.\ref{Propo2.1.v}, we obtain $| \zeta - \zeta_{n}| \leq 2|q_{2t_n+2u_n}|^{-2}$ (in this estimate we use the triangle inequality as in the proof of Equation \eqref{Eq:Teo1.1:01}). If $\lambda$ is as in Proposition \ref{Propo2.1}.\ref{Propo2.1.vii}, we can use  parts \ref{Propo2.1.ii} and \ref{Propo2.1.v} of Proposition \ref{Propo2.1} to obtain
\begin{align}
|P_n(\zeta)| 
&= |P_n(\zeta) - P_n(\zeta_{n})| \nonumber\\
&= 
\left| \zeta - \zeta_{n}\right|
\left| \begin{vmatrix} q_{2w_n-1}\zeta-p_{2w_n-1} & q_{2t_n-1} \\ q_{2w_{n}}\zeta-p_{2w_{n}} & q_{2t_n} \end{vmatrix}
+ 
\begin{vmatrix} q_{2w_n-1} & q_{2t_n-1}\zeta_{n} - p_{2t_n-1} \\ q_{2w_n} & q_{2t_n}\zeta_{ n } - p_{2t_n} \end{vmatrix} \right| \nonumber\\
&
\leq 
2\lambda \left| \zeta - \zeta_{ n }\right| \left(\left| \frac{q_{2t_n}}{q_{2w_n}}\right| 
+ 
\left|\frac{q_{2w_n}}{q_{2t_n}}\right|\right) 
\leq 
\kappa_1 \frac{|q_{2t_n}|}{|q_{2w_n}q^2_{2t_n+2u_n}|}.\label{Eq:EstimPn}
\end{align}

Let $(\bfx_n)_{n\geq 1}$, $\bfx_n=(x_{n,1}, x_{n,2}, x_{n,3}, x_{n,4})$, be the sequence in $\Za[i]^4$ given by
\[
\bfx_n:=\left(
\begin{vmatrix}
q_{2w_n - 1} & q_{2t_n - 1} \\
q_{2w_n} & q_{2t_n} \\
\end{vmatrix}, \;
\begin{vmatrix}
q_{2w_n - 1} & p_{2t_n - 1} \\
q_{2w_n} & p_{2t_n} \\
\end{vmatrix}, \;
\begin{vmatrix}
p_{2w_n - 1} & q_{2t_n - 1} \\
p_{2w_n} & q_{2t_n} \\
\end{vmatrix}, \;
\begin{vmatrix}
p_{2w_n - 1} & p_{2t_n - 1} \\
p_{2w_n} & p_{2t_n} \\
\end{vmatrix} \right).
\]
Then, $\|\bfx_n\| \leq 2|q_{2w_n} q_{2t_n}|$ for all $n\in\Na$. If $M=\sup\{B_n:n\in\Na\}$, for all $n\in\Na$ we have
\[
\left|\zeta x_{n,1} - x_{n,2}\right|
\leq 
(M + 2) \left| \frac{q_{2w_n}}{q_{2t_n}}\right|.
\]
Certainly, by Proposition \ref{Propo2.1}. \ref{Propo2.1.ii}, \ref{Propo2.1.v}  and Equation \eqref{Eq:pnqn},
\begin{align*}
    \left|\zeta x_{n,1} - x_{n,2}\right|
&= 
\begin{vmatrix}
    q_{2w_n-1} & \zeta q_{2t_n-1} - p_{2t_n-1} \\
    q_{2w_n} & \zeta q_{2t_n} - p_{2t_n} \\
\end{vmatrix}
\\
&=
\left| q_{2w_n-1}(\zeta q_{2t_n} - p_{2t_n}) - q_{2w_n}(\zeta q_{2t_n-1} - p_{2t_n-1})\right| \\
&\leq 
\left| \frac{q_{2w_n-1} }{q_{2t_n}}\right| + \left| \frac{q_{2w_n} }{q_{2t_n-1}} \right| 
\leq 
(M+2) \left| \frac{q_{2w_n} }{q_{2t_n}} \right|.
\end{align*}
We can obtain a similar bound for $|\zeta x_{n,1} - x_{n,3}|$. Consider the linear forms $\scL_{1,1}$, $\scL_{1,2}$, $\scL_{1,3}$, $\scL_{1,4}$ in the variable $\bfX=(X_1,X_2,X_3,X_4)$ and the linear forms $\scM_{1,1}$, $\scM_{1,2}$, $\scM_{1,3}$, $\scM_{1,4}$ in the variable $\bfY=(Y_1,Y_2,Y_3,Y_4)$ given by
\begin{align*}
\scL_{1,1}(\bfX) &= \zeta^2 X_1 - \zeta(X_2+X_3) + X_4, &&\quad \scM_{1,1}(\bfY) = \overline{\scL_{1,1}(\overline{ \bfY })}, \\
\scL_{1,2}(\bfX) &= \zeta X_1 - X_2, &&\quad \scM_{1,2}(\bfY) = \overline{\scL_{1,2}(\overline{ \bfY })},\\
\scL_{1,3}(\bfX) &= \zeta X_1 - X_3, &&\quad \scM_{1,3}(\bfY) = \overline{\scL_{1,3}(\overline{ \bfY })}, \\
\scL_{1,4}(\bfX) &= X_1, &&\quad \scM_{1,4}(\bfY) = \overline{\scL_{1,4}(\overline{ \bfY })}.
\end{align*}
Take $\veps>0$ as in Lemma \ref{LEM:TRANS_EPSILON}. Then, for some positive constants $\kappa_2,  \kappa_3, \kappa_4$ and all $n\in\Na$ we have
\begin{align*}
\left| \scL_{1,1}(\bfx_n) \scL_{1,2}(\bfx_n) \scL_{1,3}(\bfx_n)\scL_{1,4}(\bfx_n) \right|
&= \left| P_n(\zeta) \right| \left| \zeta x_{n,1} - x_{n,2}\right|\left| \zeta x_{n,1} - x_{n,3}\right|\left| x_{n,1}\right| 
\nonumber\\
&\leq \kappa_2 \left| \frac{q_{2t_n}}{q_{2w_n}q^2_{2t_n + 2u_n}}\frac{q_{2w_n}}{q_{2t_n}}\frac{q_{2t_n}}{q_{2w_n}}q_{2w_n}q_{2t_n} \right| \nonumber\\
&= \kappa_2 \left| \frac{q_{2t_n}^2}{q_{2t_n + 2u_n}^2}\right| \nonumber\\
&\leq \frac{\kappa_3}{\psi^{4u_n}}
\leq \frac{\kappa_3 }{|q_{2w_n}q_{2t_n}|^{\veps}} 
\leq  
\frac{\kappa_4}{\|\bfx_n\|^{\veps}}. \nonumber
\end{align*}
In the last line, the first inequality uses Proposition \ref{Propo2.1}.\ref{Propo2.1.viii}. The definition of $\scM_{1,1}$, $\scM_{1,2}$, $\scM_{1,3}$, $\scM_{1,4}$, tells us that for each $n\in\Na$ we have
\[
\left| \scM_{1,1}(\overline{\bfx_n}) \scM_{1,2}(\overline{\bfx_n}) \scM_{1,3}(\overline{\bfx_n})\scM_{1,4}(\overline{\bfx_n}) \right|
=
\left| \scL_{1,1}(\bfx_n) \scL_{1,2}(\bfx_n) \scL_{1,3}(\bfx_n)\scL_{1,4}(\bfx_n) \right|.
\]
Therefore, by Theorem \ref{TEO:SST}, there is some $\bfx=(x_1,x_2,x_3,x_4)\in\Za[i]^4$, $\bfx\neq (0,0,0,0)$, and an infinite set $\clN_1\subseteq \Na$ such that
\begin{equation}\label{Eq:Trans:01}
\forall n\in \clN_1
\quad
x_1 x_{n,1} + x_2 x_{n,2} + x_3 x_{n,3} + x_4 x_{n,4}=0.
\end{equation}
Let $(Q_n)_{n\geq 1}$ and $(R_n)_{n\geq 1}$ be the sequences given by
\[
\forall n\in\Na \qquad
Q_n:= \frac{q_{2w_n-1}q_{2t_n}}{q_{2w_n}q_{2t_n-1}}, \quad
R_n:= \zeta - \frac{p_{2n}}{q_{2n}}.
\]
For each $n\in \clN_1$, we divide Equation \eqref{Eq:Trans:01} by $q_{2w_n}q_{2t_n-1}$ to obtain
\begin{align*}
0 &= x_1(Q_n-1) + x_2\left(Q_n(\zeta - R_{2t_n}) - (\zeta - R_{2t_n-1})\right) + \nonumber\\
&\;+ x_3\left( Q_n(\zeta - R_{2w_n-1}) - (\zeta - R_{2w_n})\right)+ \nonumber\\
&\;+x_4\left( Q_n(\zeta - R_{2w_n-1})(\zeta - R_{2t_n}) - (\zeta - R_{2w_n})(\zeta - R_{2t_n-1}) \right).
\end{align*}
For instance, the coefficient of $x_2$ is 
\[
\frac{x_{n,2}}{q_{2w_n}q_{2t_n-1}}
=
Q_n\frac{p_{2t_n}}{q_{2t_n}} - \frac{p_{2t_{n}-1}}{q_{2t_n-1}}
=
Q_n(\zeta - R_{2tn}) - (\zeta - R_{2t_n-1}).
\]
Since $R_n\to 0$ as $n\to \infty$, we conclude that
\[
\lim_{\substack{n\to\infty\\ n\in \clN_1}} (Q_n-1) \left(x_1 + (x_2+x_3)\zeta + \zeta^2 x_4\right)=0.
\]
Let us show that
\begin{equation}\label{Ec:App02:proof}
x_1 + (x_2+x_3)\zeta + \zeta^2 x_4=0.
\end{equation}
By Equation \eqref{Eq:MirrorFormula}, Proposition \ref{Propo2.1}.\ref{Propo2.1.ii}, and the boundedness of $\bfB$, we conclude that
\[
0
<
\liminf_{n\to\infty} \left| \frac{q_{2w_n-1}}{q_{2w_n}}\right|
<
\limsup_{n\to\infty} \left| \frac{q_{2w_n-1}}{q_{2w_n}}\right|
<
\infty,
\]
and the same for $(q_{2t_n-1}/q_{2t_n})_{n\geq 1}$. Thus, for some infinite subset $\clN_1'\subseteq \clN_1$ and some constants $a,b\in \Za$, $a\neq b$, we have $a=B_{w_n}$ and $b=B_{t_n}$ for all $n\in \clN_1'$ and the next limits exist:
\[
\gamma_1:=
\lim_{\substack{n\to\infty \\ n\in \clN_1'}} \frac{q_{2w_n}}{q_{2w_n-1}}, 
\quad
\gamma_2:=
\lim_{\substack{n\to\infty \\ n\in \clN_1'}} \frac{q_{2t_n}}{q_{2t_n-1}}. 
\]
By \eqref{Eq:Trans:00}, $\gamma_1\in \tau_{ia}\left[ \oclC_1(-2)\cup \oclC_1(2)\right]$ and $\gamma_2\in \tau_{ib}\left[ \oclC_1(-2)\cup \oclC_1(2)\right]$ and, since these compact sets are disjoint (see Figure 
 \ref{Fig-PartmfF}), $(Q_n)_{n\in\clN_1}$ is bounded away from $1$ for large $n\in\clN_1$. Hence, \eqref{Ec:App02:proof} holds and, by $[\QU(i,\zeta),\QU(i)]\geq 3$, 
\[
x_1=x_2+x_3=x_4=0,
\;\text{ so }\;
x_2=-x_3\neq 0.
\]
As a consequence, by \eqref{Eq:Trans:01}, for $n\in\clN_1$ we have $x_{n,2}=x_{n,3}$ and
\begin{equation}\label{tredive}
P_n(X) =
\begin{vmatrix} q_{2w_n-1} & q_{2t_n-1}\\ q_{2w_n} & q_{2t_n} \end{vmatrix} X^2 - 2 \begin{vmatrix} q_{2w_n-1} & p_{2t_n-1}\\ q_{2w_n} & p_{2t_n} \end{vmatrix} X + \begin{vmatrix} p_{2w_n-1} & p_{2t_n-1}\\ p_{2w_n} & p_{2t_n} \end{vmatrix}. 
\end{equation}
Let $(\bfv_n)_{n\geq 1}$, $\bfv_n=(v_{n,1},v_{n,2},v_{n,3})\in\Za[i]^3$, be given by
\[
\forall n\in\Na\quad
\bfv_n:=\left( \begin{vmatrix} q_{2w_n-1} & q_{2t_n-1}\\ q_{2w_n} & q_{2t_n} \end{vmatrix}, \begin{vmatrix} q_{2w_n-1} & p_{2t_n-1}\\ q_{2w_n} & p_{2t_n} \end{vmatrix}, \begin{vmatrix} p_{2w_n-1} & p_{2t_n-1}\\ p_{2w_n} & p_{2t_n} \end{vmatrix}\right),
\]
so $\|\bfv_n\| \leq 2|q_{2w_n}q_{2t_n}|$. Consider the linear forms $\scL_{2,1}$, $\scL_{2,2}$, $\scL_{2,3}$ in the variable $\bfX=(X_1,X_2,X_3)$ and the linear forms $\scM_{2,1}$, $\scM_{2,2}$, $\scM_{2,3}$ in the variable $\bfY=(Y_1,Y_2,Y_3)$ given by
\begin{align*}
\scL_{2,1}(\bfX) &= \zeta^2 X_1 - 2\zeta X_2 + X_3, &&\quad \scM_{2,1}(\bfY) = \overline{\scL_{2,1}(\overline{ \bfY })}, \\
\scL_{2,2}(\bfX) &= \zeta X_1 - X_2, &&\quad \scM_{2,2}(\bfY) = \overline{\scL_{2,2}(\overline{ \bfY }),}\\
\scL_{2,3}(\bfX) &= X_1 , &&\quad \scM_{2,3}(\bfY) = \overline{\scL_{2,3}(\overline{ \bfY }).} 
\end{align*}
For all $n\in\Na$ we have
\[
\left|\scM_{2,1}(\overline{\bfv_n})\scM_{2,2}(\overline{\bfv_n})\scM_{2,3}(\overline{\bfv_n}) \right|
= 
\left|\scL_{2,1}(\bfv_n)\scL_{2,2}(\bfv_n)\scL_{2,3}(\bfv_n) \right|.
\]
By \eqref{Eq:EstimPn}, \eqref{tredive}, and parts \ref{Propo2.1.ii} and \ref{Propo2.1.v} of Proposition \ref{Propo2.1}, there are some constants $\kappa_5, \kappa_6, \kappa_7, \kappa_8>0$ depending on $\zeta$ such that for all $n\in\clN_1$ we have
\begin{align*}
\left|\scL_{2,1}(\bfv_n)\scL_{2,2}(\bfv_n)\scL_{2,3}(\bfv_n) \right|
&=
|P_n(\zeta) (\zeta v_{n,1}-v_{n,2}) v_{n,1}| \nonumber\\
&\leq 
\kappa_5 \frac{|q_{2t_n}|}{|q_{2w_n}q^2_{2t_n+2u_n}|} 
\begin{vmatrix}
q_{2w_n-1} & q_{2t_n-1}\zeta - p_{2t_n-1} \\ 
q_{2w_n}  & q_{2t_n}\zeta - p_{2t_n}
\end{vmatrix}
|q_{2w_n}q_{2t_n}| \\
&\leq 
\kappa_6 \frac{|q_{2w_n}|}{|q_{2t_n + 2u_n}|} 
\frac{|q_{2t_n}|}{|q_{2t_n + 2u_n}|} \\
&\leq  
\kappa_6 
\frac{|q_{2t_n}|}{|q_{2t_n + 2u_n}|} \leq 
\frac{\kappa_7 }{|q_{2w_n}q_{2t_n}|^{\veps}} 
\leq  
\frac{\kappa_8 }{\|\bfv_n\|^{\veps}}.  \nonumber
\end{align*}
The next to last inequality uses Lemma \ref{LEM:TRANS_EPSILON} and Proposition \ref{Propo2.1}.\ref{Propo2.1.viii}. Hence, Theorem \ref{TEO:SST} implies the existence of some infinite set $\clN_2\subseteq \clN_1$ and some $\mathbf{r}=(r_1,r_2,r_3)\in \Za[i]^3$, $\mathbf{r} \neq (0,0,0)$, such that
\[
\forall n\in \clN_2
\quad
r_1v_{n,1} + r_2v_{n,2} + r_3v_{n,3} = 0.
\]
Dividing by $q_{2w_n}q_{2t_n-1}$ we obtain
\[
\lim_{\substack{ n\to\infty \\ n\in \clN_2}} 
(Q_n-1)(r_3\zeta^2 + r_2 \zeta + r_1) = 0.
\]
Since $(Q_n)_{n\in\clN_1}$ is bounded away from $1$, we conclude
\[
r_3\zeta^2 + r_2 \zeta + r_1 = 0.
\]
This implies $[\QU(i,\zeta):\QU(i)]\leq 2$, contradicting our hypothesis. Therefore, $\zeta$ is transcendental.
\end{proof}

\section{Proof of Proposition \ref{LEM:Algorithm}}\label{SECTION:PFALGORITHM}
Recall that Proposition \ref{LEM:Algorithm} says that for each $\bfa\in\sfIr$ there is some $\bfb\in\overline{\sfR}$ such that $\overline{\Lambda}(\bfb)=\Lambda(\bfa)$. To prove it, we first need some auxiliary results. 

\subsection{Auxiliary lemmas}\label{Subsec:PrelLim}
First, we gather some technical observations. 
\begin{lem01}\label{LEM:ALGORITHM}
Take $n\in\Na$, $(x_1,\ldots,x_n)\in\scD^n$, $m \in\Za$ with $|m|\geq 2$, and $c\in\{ m, im\}$.
\begin{enumerate}[label= \roman*.]
\item If $(x_1,\ldots, x_n,c)\in\sfR(n+1)$, then $\mfFc_{n+1}(x_1,\ldots, x_n,c)=\mfFc_n(c)$.
\item For all $k\in\Na$ and $(d_1,\ldots,d_k)\in\scD^k$, if $(x_1,\ldots,x_n,c,d_1,\ldots,d_k)\in \sfR(n+k+1)$ then 
\[
\mfFc_{n + k +1}(x_1,\ldots,x_n,c,d_1,\ldots,d_k)
=
\mfFc_k(d_1,\ldots, d_k).
\]
\end{enumerate}
\end{lem01}
\begin{proof}
\begin{enumerate}[label= \roman*.]
\item The proof amounts to checking all the possible choices for $\mfFc_n(x_1,\ldots, x_n)$ and $c\in\scD$. For example, if $\mfFc_n(x_1,\ldots, x_n)=\mfFc_1(2)$, then $c=im$ for some $m\in\Za$ with $|m|\geq 2$ or $c=m$ for some $m\in\Na_{\geq 2}$ (see Figure \ref{FigB5iii}); in either case, $\mfFc_{n +1}(x_1,\ldots,x_n,c)=\mfFc(c)$. 
\item Assume $|c|\geq 3$. The previous part yields $\mfFc_{n+1}(x_1,\ldots, x_n,c)=\mfFc$ and, by \eqref{Eq:ObsDefCyls},
\begin{align}
\mfFc_{n+k+1}(x_1,\ldots, x_n,c,d_1,\ldots, d_k) 
&= T^{n+k+1}\left[ \clCc_{n + k +1}(x_1,\ldots, x_n,c,d_1,\ldots,d_k) \right] \nonumber \\
&= T^k\left[ \mfFc_{n+1}(x_1,\ldots, x_n,c) \cap \clCc_{k}(d_1,d_2,\ldots,d_k) \right] \nonumber  \\
&= T^k\left[  \mfFc \cap \clCc_{k}(d_1,\ldots,d_k) \right] \nonumber  \\
&= \mfFc_{k}(d_1,\ldots,d_k). \label{Eq:Lem7.1:Dem:00}
\end{align}
\end{enumerate}

If $|c|=2$, then $\mfFc_{n+1}(x_1,\ldots, x_n,c)\neq \vac$ may take 4 possible forms. Checking case by case, we may see that $(x_1,\ldots, x_n,c,d_1)\in\sfR(n+2)$ implies
\begin{equation}\label{Eq:Lem7.1:Dem:01}
\mfFc_{n+2}(x_1,\ldots, x_n,c, d_1) =\mfFc_1(d_1).
\end{equation}
For instance, if $c=2$, then 
\[
\iota[\mfFc_{n+1}(x_1,\ldots, x_n,2)]
=
\iota\left[\mfFc_1(2)\right]
=
\iota[\mfFc]
\cap 
\left\{z\in \Cx: \real(z) >  -\frac{1}{2} \right\}.
\]
Thus, we get \eqref{Eq:Lem7.1:Dem:01}, because for all $d_1\in\scD$ we have 
\[
\tau_{d_1}[\mfFc]\cap \iota[\mfFc_{n+1}(x_1,\ldots, x_n,c)]
=
\begin{cases}
\vac, &\text{ if } \real(d_1)< 0,\\
\tau_{d_1}[\mfFc]\cap \iota[\mfFc]\neq \vac, &\text{ if } \real(d_1)\geq 0,
\end{cases}
\]
(cfr. Figures \ref{Fig-OpenProt}, \ref{FigB5iii}). We are done if $k=1$. If not, we proceed as in \eqref{Eq:Lem7.1:Dem:00} and  using \eqref{Eq:Lem7.1:Dem:01}.
\end{proof}


\begin{lem01}\label{LEM:SEC:MIR}
If $a\in \{1+im,m+i\}$ for some $m\in\Za$, $|m|\geq 2$, then
\[
(\Mir_1\circ S)(a)=(S\circ\Mir_1)(a)
\quad\text{ and } \quad
(\Mir_2\circ S)(a)=(S\circ\Mir_2)(a).
\] 
\end{lem01}
\begin{proof}
The lemma follows from the definitions of $\Mir_1$, $\Mir_2$, and $S$.
\end{proof}

\begin{lem01}\label{LEM:ALGORITHM03}
Take $\bfb=\sabu\in \Omega$. If $j\in\Na$ satisfies $(b_1,\ldots, b_j)\in \sfR(j)$ and $(b_1,\ldots, b_j)\in \sfIr(j+1)$, then 
\[
\real([b_j;b_{j+1},b_{j+1},\ldots]_{\Cx})= \frac{1}{2}
\quad\text{ or }\quad
\imag([b_j;b_{j+1},b_{j+1},\ldots]_{\Cx})= \frac{1}{2}.
\]
\end{lem01}
\begin{proof}
This is a consequence of Lemma \ref{Lema:LemaB5} and its proof.
\end{proof}

\subsection{Proof of Proposition \ref{LEM:Algorithm}}\label{SUBSECTION:PFPROC}
We deduce Proposition \ref{LEM:Algorithm} after showing that Procedure \ref{Algor} gives the desired output. To this end, let $\bfa\in\sfIr$ be arbitrary, and $(\bfb_N)_{N\geq 1}$ and $\bfb$ as in Procedure \ref{Algor}. We will verify the next \textit{three assertions}:
\begin{enumerate}[label= \Alph*., ref=\Alph*]
\item \label{AsserA} The sequence $(\bfb_N)_{N\geq 0}$ converges in the product topology to some $\bfb\in \scD^{\Na}$,
\item \label{AsserB} $\bfb:=\displaystyle\lim_{N\to \infty}\bfb_N\in\overline{\sfR}$, and 
\item \label{AsserC} $\overline{\Lambda}(\bfb)=\Lambda(\bfa)$.
\end{enumerate}

First, we show that the sequences $\bfb_{N}$ and $\bfb_{N+1}$ have increasingly longer regular prefixes in common; this will be proved below by induction. With that, Assertions \ref{AsserA} and \ref{AsserB} will follow. Our arguments will also imply that
\[
\Lambda(\bfa)=[0;b_N(1),b_N(2),b_N(3),\ldots],
\] 
which will lead us to conclude Assertion \ref{AsserC}.

\subsubsection*{Base of induction.} Let us start with $N=0$. Since every cylinder of level $1$ is regular, there is some $j_0\in\Na$ such that
\[
(a_1,\ldots,a_{j_0})\in\sfR(j_0) 
\text{ and }
(a_1,\ldots,a_{j_0}, a_{j_0+1})\not\in\sfR(j_0+1). 
\]
Note that $\mfFc_{j_0}(a_1,\ldots,a_{j_0})\neq \mfFc$, otherwise $(a_1,\ldots,a_{j_0}, a_{j_0+1})$ would be regular. 
By the proof of Lemma \ref{Lema:LemaB5}, we may conclude the existence of some $m\in\Za$ such that 
\[
|m|\geq 2
\quad\text{ and }\quad
a_{j_0+1} \in\left\{ 1+im, m+i\right\}.
\]
Put $t=1$ if $a_{j_0+1}=m+i$ and $t=2$ if $a_{j_0+1}=1+m$. We define 
\[
\bfb_1
:=
(a_1,\ldots, a_{j_0}, S(a_{j_0 +1 }), \Mir_t(a_{j_0 + 2}), \Mir_t(a_{j_0 + 3}), \Mir_t(a_{j_0 + 4}), \ldots ).
\]
We verify that the sequence $\bfb_1=(b_1(j))_{j\geq 1}$ has the next properties:
\begin{align}
(b_1(1),\ldots,  b_1(j_0+1)) \in \sfR(j_0+1), \label{Eq:Dem:Algor:Prob:b1:01} \\
\Lambda(\bfa) = [0;b_1(1),b_1(2),b_1(3),\ldots]. \label{Eq:Dem:Algor:Prob:b1:02} 
 \end{align}

{\bf Case $|m|\geq 3$.} By Lemma \ref{Lema:LemaB5}, $\mfFc_{j_0+1}\left(a_1,\ldots,a_{j_0}, S(a_{j_0+1})\right) = \mfFc $ and \eqref{Eq:Dem:Algor:Prob:b1:01} follows. Now we show \eqref{Eq:Dem:Algor:Prob:b1:02}. If $a_{j_0+1}=1 + im$, then Lemma \ref{LEM:ALGORITHM03} and $[0;a_{j_0+2},a_{j_0+3},\ldots]_{\Cx}\in \mfF$ give
\begin{equation}\label{Eq:AlgR12N0}
\real\left( [a_{j_0 + 1};a_{j_0 + 2},a_{j_0 + 3},\ldots]_{\Cx}\right)
=
\frac{1}{2}.
\end{equation}
Since $z= 1 - \overline{z}$ holds whenever $\real(z)=\frac{1}{2}$, we obtain 
\begin{align*}
[a_{j_0 + 1};a_{j_0 + 2},a_{j_0 + 3},\ldots]_{\Cx}
&=
1 - \overline{[a_{j_0 + 1};a_{j_0 + 2},a_{j_0 + 3},\ldots]_{\Cx}} \\
&=
[ 1 - \overline{a_{j_0 + 1}}; -\overline{a_{j_0 + 2}}, -\overline{a_{j_0 + 3}},\ldots ] \\
&=
[ S(a_{j_0+1}) ; \Mir_2(a_{j_0 + 1}) , \Mir_2(a_{j_0 + 2}), \Mir_2(a_{j_0 + 3}) , \ldots ].
\end{align*}
Then, we may conclude \eqref{Eq:Dem:Algor:Prob:b1:02}:
\begin{align*}
\Lambda(\bfa)
&=
[0 ;a_1, \ldots, a_{j_0}, a_{j_0+1}, a_{j_0+2}, \ldots]_{\Cx} \\
&=
[ 0;a_1,a_2,\ldots, a_{j_0}, S(a_{j_0+1}), \Mir_2(a_{j_0+2}), \Mir_2(a_{j_0+3})\ldots ] \\
&=
[ 0;b_1(1), b_1(2),\ldots, b_1(j_0), b_1(j_0+1), b_1(j_0+2), b_1(j_0+3),\ldots].
\end{align*}
Similarly, if $a_{j_0+1}=m + i$, from Lemma \ref{LEM:ALGORITHM03} and $[0;a_{j_0+2},a_{j_0+3},\ldots]_{\Cx}\in \mfF$ we get
\begin{equation}\label{Eq:AlgI12N0}
\imag\left( [a_{j_0 + 1};a_{j_0 + 2},a_{j_0 + 3},\ldots]_{\Cx}\right)
=
\frac{1}{2}.
\end{equation}
Since $z= \overline{z} +i$ when $\imag(z)=\frac{1}{2}$, we have
\begin{align*}
[a_{j_0 + 1};a_{j_0 + 2},a_{j_0 + 3},\ldots]_{\Cx}
&=
i + \overline{[a_{j_0 + 1};a_{j_0 + 2},a_{j_0 + 3},\ldots]_{\Cx}} \\
&=
[ i + \overline{a_{j_0 + 1}} ;\overline{a_{j_0 + 2}}, \overline{a_{j_0 + 3}},\ldots ] \\
&=
[ S\left( a_{j_0 + 1}\right); \Mir_1(a_{j_0 + 1}), \Mir_1(a_{j_0 + 2}), \Mir_1(a_{j_0 + 3}) , \ldots]
\end{align*}
and \eqref{Eq:Dem:Algor:Prob:b1:01} follows.

{\bf Case $|m|= 2$.} Apply $\iota$ on the sets in Lemma \ref{Lema:LemaB5} to obtain $a_{j_0+1}\in\{1+2i,1-2i, 2+i, -2+i\}$. Suppose that $a_{j_0+1}=1+2i$, then either
\[
\mfFc_{j_0}(a_1,\ldots, a_{j_0}) = \mfFc_1(-1-i)
\;\text{ or }\;
\mfFc_{j_0}(a_1,\ldots, a_{j_0}) = \mfFc_1(-2)
\]
(the only sets in Lemma \ref{Lema:LemaB5} leading to an irregular cylinder when the next partial quotient is $1+2i$ are $\mfFc_1(-1-i)$ and $\mfFc_1(-2)$). In either case, we have 
\[
\mfFc_{j_0+1}\left( a_1,\ldots, a_{j_0}, S(a_{j_0+1})\right)
=
\mfFc_{j_0+1}(a_1,\ldots, a_{j_0}, 2i)
=
\mfFc_1(2i)\neq\vac
\]
and we get \eqref{Eq:Dem:Algor:Prob:b1:01}. We show \eqref{Eq:Dem:Algor:Prob:b1:02} as in the case $|m|\geq 3$. The possibilities $a_{j_0+1}\in\{1-2i, 2+i,2-i\}$ are treated similarly.

\subsubsection*{Inductive step} Assume that for some $N\in\Na$ we can perform $N$ iterations of the procedure and that for each $k\in\{1,\ldots,N\}$ we have
\begin{equation}\label{Eq:IH-Reg} 
(b_k(1),\ldots, b_k(j_{k-1}+1)) \in \sfR(j_{k-1}+1),
\end{equation}
\begin{equation}\label{Eq:IH-ReIm}
\real\left([a_{j_{k-1} + 1}; a_{j_{k-1} + 2}, \ldots ]_{\Cx}\right) = \frac{1}{2} 
 \;\text{ or }\;
\imag\left([a_{j_{k-1} + 1}; a_{j_{k-1} + 2} , \ldots ]_{\Cx}\right) = \frac{1}{2} \text{ and },
\end{equation}
\begin{equation}\label{Eq:IH-cf} 
\Lambda(\bfa)
 =
[0;b_k(1),b_k(2),b_k(3),\ldots].
\end{equation}
\begin{rema}\label{Remark:01}
    By the definition of $\bfb_N$, for each $M\in\{1,\ldots, N\}$ there exists some $\Mir\in\MIR$ such that 
    \[
    \forall n\in\Na \quad
    b_M(j_{M-1} + 1 + n) = \Mir\left( b_{M-1}(j_{M-1} + 1 + n)\right).
    \]
    As a consequence, there is some $\Mir'\in\MIR$ for which
    \[
    \forall n\in\Na \quad
    b_M(j_{M-1} + 1 + n) = \Mir'\left( a_{ j_{M-1} + 1 + n }\right).
    \]
\end{rema}
 Let us discuss the $(N+1)$th iteration. Suppose there is some $j_{N}\in\Na$ satisfying
\begin{equation}\label{Eq:InductiveStep:01}
\left( b_N(1),\ldots, b_N(j_N)\right)\in \sfR(j_{N})
\;\text{ and }\;
\left( b_N(1),\ldots, b_N(j_N+1)\right)\not\in \sfR(j_{N}+1).    
\end{equation}
Assertions \ref{AsserA}, \ref{AsserB}, and \ref{AsserC} are trivially true if no such $j_N$ exists. We claim that 
\begin{equation}\label{Eq:InductiveStep:02}
    \left( b_N(1),\ldots, b_N(j_N),  S(b_N(j_N+1)) \right) \in \sfR(j_{N}+1), 
\end{equation}
\begin{equation}\label{Eq:InductiveStep:03}
    \real\left([a_{j_{N} + 1};a_{j_{N} + 2} , \ldots ]_{\Cx}\right) = \frac{1}{2}
\quad \text{ or } \quad
\imag\left([a_{j_{N} + 1};a_{j_{N} + 2},\ldots ]_{\Cx}\right) = \frac{1}{2}, \text{ and}
\end{equation}
\begin{equation}\label{Eq:InductiveStep:03.5}
    \Lambda(\bfa) 
    = 
    [0;b_{N+1}(1),b_{N+1}(2),b_{N+1}(3),\ldots].
\end{equation}
From \eqref{Eq:IH-Reg}, we deduce $j_N-j_{N-1}\geq 1$. We consider two cases: $j_N-j_{N-1}= 1$ and $j_N-j_{N-1}\geq 2$.


\textbf{Case $j_{N}=j_{N-1} + 1$.}  Since $b_N(j_N) = S(b_{N-1}(j_{N}))= S(b_{N-1}(j_{N-1}+1))$, there is some $m\in \Za$, $|m|\geq 2$, such that $b_N(j_N)\in \{m, im,-m,-im\}$. If we had $|m|\geq 3$, then Lemma \ref{LEM:ALGORITHM} would give
\[
\mfFc_{j_{N}}(b_N(1), \ldots, b_N(j_N))
=
\mfFc_{j_{N-1}+1}(b_{N-1}(1), \ldots, S(b_{N-1}(j_{N-1}+1))
=
\mfFc,
\]
which implies $(b_N(1), \ldots, b_N(j_N),b_N(j_N+1))\in \sfR(j_N+1)$, contradicting \eqref{Eq:InductiveStep:01}. Therefore, $|m|=2$ and $b_N(j_N)\in \{-2,-2i,2,2i\}$. Let us further assume that $b_N(j_N)=-2$, the other possibilities being treated similarly. Then, we have 
\[
b_{N-1}(j_{N})
=
b_{N-1}(j_{N-1}+1)
\in
\{-2 + i, -2 - i\}.
\]
Suppose that $b_{N-1}(j_{N})=-2+i$, an analogous argument works for $-2-i$. Then, Condition \eqref{Eq:InductiveStep:01} and Lemma \ref{LEM:ALGORITHM} imply
\begin{equation}\label{Eq:InductiveStep:04}
\mfFc_{j_N}\left(b_N(1), \ldots, b_N(j_N)\right)= \mfFc_1(-2)
\quad\text{ and }\quad
\real\left(b_N(j_N+1)\right)\geq 1 .
\end{equation}
To see the inequality in \eqref{Eq:InductiveStep:04}, note that $\mfFc_1(-2)=-\mfFc_1(2)$, $\iota[\mfFc_1(-2)]=-\iota[\mfF(2)]$ and see Figure \ref{FigB5iii}. We claim that 
\begin{equation}\label{Eq:InductiveStep:05}
    \exists \Mir\in\MIR\quad
    \forall n\in\Na \quad
    b_{N-1}(j_{N-1}+n) = \Mir(a_{j_{N-1}+n }).
\end{equation}
Certainly, when $N=1$, we consider $\Mir$ to be the identity map on $\Cx$. For $N>1$, the induction hypothesis gives $j_{N-1} \geq j_{N-2}+1$, so $j_{N-1}+n \geq j_{N-2}+n+1$ for all $n\in\Na$ and \eqref{Eq:InductiveStep:05} follows from Remark \ref{Remark:01}. 

Next, we separately consider every possibility for $\Mir$ in \eqref{Eq:InductiveStep:05}. Recall that
\[
-2+i = b_{N-1}(j_N)=\Mir(a_{j_N}).
\]

\begin{figure}[ht!]
\begin{center}
\includegraphics[scale=0.85,  trim={6cm 19.25cm 3.5cm 5.75cm},clip]{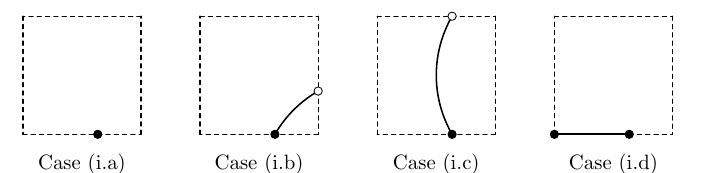}
\caption{ Possible forms of $\mfF(a_1,\ldots,a_{j_N})$ when $\Mir$ is the identity map. \label{Fig-Alg02}}
\includegraphics[scale=0.85,  trim={6cm 19.25cm 3.5cm 5.75cm},clip]{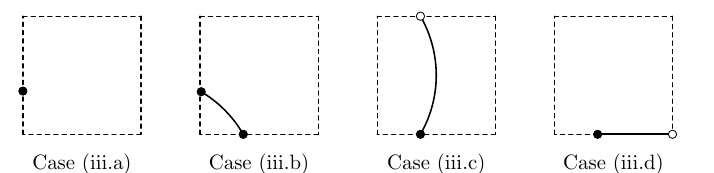}
\caption{ Possible forms of $\mfF(a_1,\ldots,a_{j_N})$ when $\Mir=\Mir_2$. \label{Fig-Alg03}}
 \end{center}
\end{figure}

\begin{enumerate}[label= \roman*.]
\item When $\Mir$ is the identity map, we have $a_{j_N}=-2+i$ and $\real(a_{j_N+1})\geq 1$. Hence, the set $\mfF_{j_N}(a_1,\ldots, a_{j_N})\in \sfIr(j_N)$ is of one of the forms depicted in Figure \ref{Fig-Alg02}. To check this observation, we determine inductively on $n$ the possible forms of $\mfF_{n}(\bfc)$, $\bfc\in\sfIr(n)$. We consider each form independently.

\begin{enumerate}[label= (i.{\alph*})., ref=(i.{\alph*}) ]
\item \label{Mir.Ident.(a)} 
Assume that $\mfF_{j_N}(a_1,\ldots, a_{j_N})=\{[0;a_{j_N+1}, a_{j_N+2},\ldots]_{\Cx}\}= \{\zeta_4\}$. This happens, for example, when $(a_1,\ldots, a_{j_N})=(2i,-2+i,1+2i,-2+i)$.
We conclude \eqref{Eq:InductiveStep:03} from $\real(\zeta_4^{-1})=\frac{1}{2}$. Moreover, we have
\begin{align*}
\left((b_{N-1}(j_N+k)\right)_{k\geq 0} &= 
(-2+i,a_{j_N+1}, a_{j_N+2}, a_{j_N+3},\ldots) \\
&=
(-2+i,1+2i,-2+i,1+2i,-2+i,\ldots). 
\end{align*}
Then, the definition of $\bfb_N$ and $j_{N}=j_{N-1}+1$ yield
\[
\left((b_{N}(j_N+k)\right)_{k\geq 0} 
=
(-2,1-2i,-2-i,1-2i,-2-i,1-2i\ldots). 
\]
Then, $j_{N+1}=j_{N}+1$ by $(-2,1-2i)\in\sfIr(2)$ and 
\[
\left((b_{N+1}(j_N+k)\right)_{k\geq 0} 
=
(-2,-2i,2-i,-1-2i,2-i,-1-2i\ldots). 
\]
We now argue as in the case $N=0$ to conclude
\[
[b_N(j_N); b_N(j_N+1) ,\ldots]
=
[b_{N+1}(j_N); b_{N+1}(j_N+1),  \ldots ].
\]
Since $b_{N}(k)=b_{N+1}(k)$ for all $k\in\{1,\ldots, j_N-1\}$, Equation \eqref{Eq:InductiveStep:03.5} follows.

\item \label{Mir.Ident.(b)} Suppose that $\mfF_{j_N}(a_1,\ldots, a_{j_N})$ is the arch of $C(1-i)$ going from $-\zeta_1=\frac{1}{2} - i\alpha$ (excluded) to $\zeta_4$ (included). An example of this case is $(a_1,\ldots, a_{j_N})=(-2,1+3i,-2+i)$. Then,  
\begin{align*}
    \left\{ w\in \Cx: \imag(w)=\frac{1}{2}\right\} \cap \iota[\mfF_{j_N}(a_1,\ldots, a_{j_N})] &= \vac,\\
    \left\{ w\in \Cx: \real(w)=\frac{1}{2}\right\} \cap \iota[\mfF_{j_N}(a_1,\ldots, a_{j_N})] &= \{\tau_{1+2i}(\zeta_2)\},
\end{align*}
so $\zeta_2=[0;a_{j_{N+1}},a_{j_{N+2}},\ldots]$ and we argue as in Case \ref{Mir.Ident.(a)}.

\item \label{Mir.Ident.(c)} When $\mfF_{j_N}(a_1,\ldots, a_{j_N})$ is the arch of $C(1)$ from $\overline{\zeta_4}$ (excluded) to $\zeta_4$ (included), we argue as in Case \ref{Mir.Ident.(b)}.
An instance of this case is $(a_1,\ldots, a_{j_N})=(2i,-3+i, 1+2i, -2+i)$.

\item \label{Mir.Ident.(d)} 
When $\mfF_{j_N}(a_1,\ldots, a_{j_N})=[-\frac{1+i}{2},\zeta_4]$, we apply $\iota$ on $\mfF_{j_N}(a_1,\ldots, a_{j_N})$ to get $a_{j_{N}+1}\in \{-1+i, -1 + 2i, 2i, 1+2i\}$ and, since $\real(a_{j_{N}+1} )\geq 1$, we must have $a_{j_{N}+1}=1+2i$. Then,
\[
\iota\left[ \mfF_{j_N}(a_1,\ldots, a_{j_N})\right]
\cap 
\tau_{a_{j_{N}+1}}[\mfF]
= 
\left\{ \tau_{a_{j_N}+1}(\zeta_2) \right\}
\]
and $\zeta_2=[0;a_{j_N+1}, a_{j_N+2},\ldots]_{\Cx}$. Now, we proceed as in Case \ref{Mir.Ident.(a)}.
An instance of this case is $(a_1,\ldots, a_{j_N})=(2i,-2+i)$.
\end{enumerate}

As we will see, the other options for $\Mir$ are impossible.
\item If $\Mir=\Mir_1$, then $a_{j_N}=-2-i$ but \eqref{Eq:IH-ReIm} fails for $k=N$.

\item When $\Mir=\Mir_2$, we argue as when $\Mir$ is the identity map so we omit some details. We have $a_{j_N}=2+i$, $\real(a_{j_N+1})\leq -1$ and, by \eqref{Eq:IH-ReIm} for $k=N$, 
\begin{equation}\label{Eq:Mir:Mir02}
    \imag([a_{j_N};a_{j_N+1},a_{j_N+2},\ldots]_{\Cx}) = \frac{1}{2}.
\end{equation}
The possible options for $\mfFc_{j_N}(a_1,\ldots, a_{j_N})$ are depicted in Figure \ref{Fig-Alg03}.
\begin{enumerate}[label= (iii.{\alph*})., ref = (iii.{\alph*}) ]
\item \label{Mir.Mir2.(a)}  If $\mfF_{j_N}(a_1,\ldots,a_{j_N})=\{\zeta_2\}=\{[0;a_{j_N+1}, a_{j_N+2},\ldots]_{\Cx}\}$, then $[a_{j_N};a_{j_N+1},\ldots]_{\Cx} = -2+i + \zeta_2$, contradicting \eqref{Eq:Mir:Mir02}.

\item \label{Mir.Mir2.(b)} If $\mfF_{j_N}(a_1,\ldots,a_{j_N})$ were the arch of $C(-1-i)$ going from $\zeta_3$ (included) to $\zeta_2$ (included), by \eqref{Eq:IH-ReIm} we would obtain $\zeta_3=[0;a_{j_N+1},a_{j_N+2}, \ldots]_{\Cx}$ so $a_{j_{N}+1}=2i$, but this contradicts $\real(a_{j_N+1})\leq -1$.
\item \label{Mir.Mir2.(c)} The impossibility of $\mfF_{j_N}(a_1,\ldots, a_{j_N})$ being the arch of $C(-1)$ from $\zeta_2$ (included) to $\overline{\zeta_2}$ (excluded) is shown as in Case \ref{Mir.Mir2.(b)}.

\item \label{Mir.Mir2.(d)} We cannot have $\mfF_{j_N}(a_1,\ldots,a_{j_N})=[\zeta_3,\frac{1-i}{2})$, because it would imply $a_{j_N+1}\in \{2i, 1+2i, 1+i\}$ which is incompatible with $\real(a_{j_N+1})\leq -1$.

\end{enumerate}

\item When $\Mir=\Mir_1\Mir_2$, then $a_{j_N}=2-i$ and $\real(a_{{j_N}+1})\leq -1$. However, we cannot have \eqref{Eq:IH-ReIm} for $k=N$ if $a_{j_N}=2-i$.

\end{enumerate}

\textbf{Case $j_{N}\geq j_{N-1} + 2$.} By the definition of $\bfb_N$, for some $m\in\Za$ with $|m|\geq 2$, we have 
\[
b_N(j_{N-1} +1 )
=
S\left( b_{N-1}(j_{N-1} +1 )\right)
\in
\left\{ m, im\right\}.
\]
Then, by Lemma \ref{LEM:ALGORITHM},
\begin{align*}
\mfFc_{j_{N}} \left( b_{N}(1), \ldots,   b_{N} (j_{N})\right)
&=
\mfFc_{j_{N} - j_{N-1} - 1 }(b_{N}(j_{N-1}+2), \ldots,  b_{N-1} (j_{N}))\neq\vac, \\
\mfFc_{j_{N}+1} \left( b_{N}(1), \ldots,  b_{N} (j_{N}+1)\right)
&=
\mfFc_{j_{N} - j_{N-1} }\left(b_{N}(j_{N-1}+2), \ldots,  b_{N} (j_{N}), b_{N} (j_{N}+1) \right) = \vac.
\end{align*}
By Remark \ref{Remark:01}, there is some $\Mir\in\MIR$ such that 
\[
\forall n\in\Na
\quad
b_N(j_{N-1} + 1 +n) 
=
\Mir\left(a_{j_{N-1} + 1 +n}\right).
\]
By the choice of $j_N$, Proposition \ref{Propo:Symm} and Lemma \ref{LEM:ALGORITHM}, we conclude
\[
\left(a_{j_{N-1} + 2}, \ldots, a_{j_{N}}\right)\in \sfR(j_N - j_{N-1}-1)
\;\text{ and }\;
\left(a_{j_{N-1} + 2}, \ldots, a_{j_{N}+1}\right)\in \sfIr(j_N - j_{N-1}-1).
\]
Applying the case $N=0$ we get
\[
\left(a_{j_{N-1} + 2}, \ldots, a_{j_{N}}, S(a_{j_{N}+1})\right)\in \sfR(j_N - j_{N-1}-1)
\]
and \eqref{Eq:InductiveStep:02} follows. We now use Lemma \ref{LEM:ALGORITHM03} to conclude \eqref{Eq:InductiveStep:03}.

For $\bfx=\saxu\in \scD^{\Na}$, put $\Mir(\bfx)=(\Mir(x_j))_{j\geq 1}$. Let $\bfd=(d_j)_{j\geq 1}$ be obtained from the first iteration of Procedure \ref{Algor} applied on $\bfa':=(a_{n})_{n\geq j_{N-1} + 2}$, then $\Lambda(\bfa')=[0;d_1,d_2,\ldots]$.  If $\bfb':=(b_{N}(n))_{n\geq j_{N-1} + 2}$ and $\bfe:=
(b_{N+1}(n))_{n\geq j_{N-1} + 2}$, then $\Mir(\bfa') = \bfb'$, $\Mir(\bfd) = \bfe$ (see Lemma \ref{LEM:SEC:MIR}) and direct computations give
\begin{align*}
\overline{\Lambda}(\Mir(\bfa'))
&=
\Mir\left(\overline{\Lambda}(\bfa')\right) \\
&=
\Mir\left([0;d_1,d_2, \ldots] \right) 
=
[0;\Mir(d_1),\Mir(d_2), \ldots]
=
[0;e_1,e_2,\ldots].
\end{align*}
Equation \eqref{Eq:InductiveStep:03.5} follows. We can, thus, conclude Assertions \ref{AsserA} and \ref{AsserB}.

\subsubsection*{Conclusion of the proof.} Lastly, we show Assertion \ref{AsserC}. Let $\sepn$, $\seqn$ be the sequences defined as in \eqref{Eq:pnqn} corresponding to $\bfb$. Then, $\frac{p_n}{q_n}\to \overline{\Lambda}(\bfb)$ as $n\to\infty$ and $(|q_n|)_{n\geq 0}$ is strictly increasing. Take $N\in\Na$. By Remark \ref{Remark:01}, there is some $\Mir'\in\MIR$ such that
\[
\xi_{N+1}
:=
[0;b_{N+1}(j_N+2),b_{N+1}(j_N+3),\ldots]
=
\Mir'\left( [0;a_{j_N+2}, a_{j_N+3},\ldots]_{\Cx}\right)
\in\omfF. 
\]
Since $(b_1,\ldots, b_{j_N+1})=(b_{N+1}(1),\ldots,b_{N+1}(j_{N}+1))$, the above argument and the elementary theory of continued fractions \cite[Theorems 2 and 5]{Khi1997} yield
\[
\Lambda(\bfa)
= 
[0;b_{N+1}(1),b_{N+1}(2),\ldots] 
=
\frac{p_{j_N+1}}{q_{j_N+1}} + \frac{(-1)^{j_N+1}\xi_{N+1}}{q_{j_{N}+1} \left(1 + \frac{q_{j_N}}{q_{j_{N}+1} } \xi_{N+1}\right)}.
\]
Then, using $|\xi_{N+1}|\leq \sqrt{2}$, we have
\[
\left|\Lambda(\bfa) - \frac{p_{j_N+1}}{q_{j_N+1}}\right|
\leq
\frac{2}{2-\sqrt{2}} \, \frac{1}{|q_{j_{N}+1}|^2}
\to 0
\;\text{ as }\; N\to\infty,
\]
and $\overline{\Lambda}(\bfb)=\Lambda(\bfa)$ follows. With the three assertions proven, Proposition \ref{LEM:Algorithm} is established.

\section{Questions for further studies}

In recent years, the theory of continued fractions in certain imaginary quadratic number fields has attracted some attention \cite{EiItoNak2019, MR4430094, EiNakadaNatsui2023}. Although our work relies on the symmetries of the HCF process, we believe that a general procedure can be obtained to cover the continued fraction algorithms used in the papers mentioned above. Hence, one of the open questions is to explore:
\begin{question}
How can the results of this paper be extended to the continued fraction algorithms in a certain number field such as those studied in \cite{EiItoNak2019, MR4430094, EiNakadaNatsui2023}?
\end{question}
Thermodynamic formalism studies measures of maximal entropy and equilibrium measures of dynamical systems. This area has interesting interpretations in physics and (fractal) dimension theory (like Bowen's formula \cite{bowen1979hausdorff}). Thermodynamic formalism for non-finite alphabets is an active area with interesting connections to surface dynamics (see the survey \cite{sarig1999thermodynamic}). Unlike the finite-alphabet theory, techniques and results depend on the class of systems. 
\begin{question}
    Are there natural equilibrium measures for  $(\overline{\mathsf{R}}, \sigma)$ and does the variational principle hold?
\end{question}

Another way to represent freedom in a dynamical system is through the space of invariant measures. The richest possible simplex of invariant measures one can attain is the Poulsen simplex. The Poulsen simplex is achieved for finite-alphabet shifts of finite type and some countable alphabet subshifts \cite{iommi2021space}. 

\begin{question}
    What can we say about the simplex of invariant measures of $(\overline{\mathsf{R}}, \sigma)$?
    \end{question}

\textbf{Acknowledgements.} 
The authors thank anonymous referees for the careful reading of the manuscript and several comments that have improved the quality of this work. We also thank Nikita Shulga and Ben Ward for their comments and suggestions. 

Research by García-Ramos was funded by the Strategic Excellence Initiative program of the Jagiellonian University with grant U1U/W16/NO/01.03. Research by González Robert was partially funded by CONAHCyT through the program \textit{Estancias posdoctorales por México}. Both González Robert and Hussain were funded by the ARC Discovery Project 200100994. 

\bibliographystyle{abbrv}
\bibliography{refHCF}

\begin{thebibliography}{10}

\bibitem{AdaBug2005}
B.~Adamczewski and Y.~Bugeaud.
\newblock On the complexity of algebraic numbers. {II}. {C}ontinued fractions.
\newblock {\em Acta Math.}, 195:1--20, 2005.

\bibitem{AdaBug2007}
B.~Adamczewski and Y.~Bugeaud.
\newblock On the complexity of algebraic numbers. {I}. {E}xpansions in integer
  bases.
\newblock {\em Ann. of Math. (2)}, 165(2):547--565, 2007.

\bibitem{adler1981construction}
R.~Adler, M.~Keane, and M.~Smorodinsky.
\newblock A construction of a normal number for the continued fraction
  transformation.
\newblock {\em J. Number Theory}, 13(1):95--105, 1981.

\bibitem{AirJacKwiMan2020}
D.~Airey, S.~Jackson, D.~Kwietniak, and B.~Mance.
\newblock Borel complexity of sets of normal numbers via generic points in
  subshifts with specification.
\newblock {\em Trans. Amer. Math. Soc.}, 373(7):4561--4584, 2020.

\bibitem{becher2019absolutely}
V.~Becher and S.~A. Yuhjtman.
\newblock On absolutely normal and continued fraction normal numbers.
\newblock {\em Int. Math. Res. Not. IMRN}, 2019(19):6136--6161, 2019.

\bibitem{bowen1971periodic}
R.~Bowen.
\newblock Periodic points and measures for {A}xiom {A} diffeomorphisms.
\newblock {\em Trans. Amer. Math. Soc.}, 154, 1971.

\bibitem{bowen1979hausdorff}
R.~Bowen.
\newblock Hausdorff dimension of quasicircles.
\newblock {\em Inst. Hautes \'Etudes Sci. Publ. Math.}, (50):11--25, 1979.

\bibitem{Bug2013}
Y.~Bugeaud.
\newblock Automatic continued fractions are transcendental or quadratic.
\newblock {\em Ann. Sci. \'Ec. Norm. Sup\'er. (4)}, 46(6):1005--1022, 2013.

\bibitem{BugGeroHus2023}
Y.~Bugeaud, G.~Gonz\'{a}lez~Robert, and M.~Hussain.
\newblock Metrical properties of {H}urwitz continued fractions.
\newblock {\em Preprint: arXiv:2306.08254}, 2023.

\bibitem{BugKim2019}
Y.~Bugeaud and D.~H. Kim.
\newblock A new complexity function, repetitions in {S}turmian words, and
  irrationality exponents of {S}turmian numbers.
\newblock {\em Trans. Amer. Math. Soc.}, 371(5):3281--3308, 2019.

\bibitem{CliTho2021}
V.~Climenhaga and D.~J. Thompson.
\newblock Beyond {B}owen's specification property.
\newblock In {\em Thermodynamic formalism}, volume 2290 of {\em Lecture Notes
  in Mathematics}, pages 3--82. Springer, Cham, [2021] \copyright 2021.

\bibitem{Dani2015}
S.~G. Dani.
\newblock Continued fraction expansions for complex numbers---a general
  approach.
\newblock {\em Acta Arith.}, 171(4):355--369, 2015.

\bibitem{DanNog2014}
S.~G. Dani and A.~Nogueira.
\newblock Continued fractions for complex numbers and values of binary
  quadratic forms.
\newblock {\em Trans. Amer. Math. Soc.}, 366(7):3553--3583, 2014.

\bibitem{deka2022borel}
K.~Deka, S.~Jackson, D.~Kwietniak, and B.~Mance.
\newblock Borel complexity of sets of points with prescribed {B}irkhoff
  averages in polish dynamical systems with a specification property.
\newblock {\em arXiv:2210.08937}, 2022.

\bibitem{DiamondVaaler1986}
H.~G. Diamond and J.~D. Vaaler.
\newblock Estimates for partial sums of continued fraction partial quotients.
\newblock {\em Pacific J. Math}, 122(1):73--82, 1986.

\bibitem{EiItoNak2019}
H.~Ei, S.~Ito, H.~Nakada, and R.~Natsui.
\newblock On the construction of the natural extension of the {H}urwitz complex
  continued fraction map.
\newblock {\em Monatsh. Math.}, 188(1):37--86, 2019.

\bibitem{MR4430094}
H.~Ei, H.~Nakada, and R.~Natsui.
\newblock On the existence of the {L}egendre constants for some complex
  continued fraction expansions over imaginary quadratic fields.
\newblock {\em J. Number Theory}, 238:106--132, 2022.

\bibitem{EiNakadaNatsui2023}
H.~Ei, H.~Nakada, and R.~Natsui.
\newblock On the ergodic theory of maps associated with the nearest integer
  complex continued fractions over imaginary quadratic fields.
\newblock {\em Discrete Contin. Dyn. Syst.}, 43(11):3883--3924, 2023.

\bibitem{Gero2020PP}
G.~Gonz\'{a}lez~Robert.
\newblock Purely periodic and transcendental complex continued fractions.
\newblock {\em Acta Arith.}, 194(3):241--265, 2020.

\bibitem{HeXio2021-01}
Y.~He and Y.~Xiong.
\newblock Sets of exact approximation order by complex rational numbers.
\newblock {\em Math. Z.}, 301:199–223, 2022.

\bibitem{Hur1887}
A.~Hurwitz.
\newblock \"{U}ber die {E}ntwicklung complexer {G}r\"{o}ssen in
  {K}ettenbr\"{u}che.
\newblock {\em Acta Math.}, 11(1-4):187--200, 1887.

\bibitem{iommi2021space}
G.~Iommi and A.~Velozo.
\newblock The space of invariant measures for countable {M}arkov shifts.
\newblock {\em J. Anal. Math.}, 143(2):461--501, 2021.

\bibitem{Kal2021}
O.~Kallenberg.
\newblock {\em Foundations of modern probability}, volume~99 of {\em
  Probability Theory and Stochastic Modelling}.
\newblock Springer, Cham, third edition, [2021] \copyright 2021.

\bibitem{Karpenkov2022}
O.~N. Karpenkov.
\newblock {\em Geometry of continued fractions}, volume~26 of {\em Algorithms
  and Computation in Mathematics}.
\newblock Springer, Berlin, second edition, [2022] \copyright 2022.

\bibitem{Kaufman1980}
R.~Kaufman.
\newblock Continued fractions and {F}ourier transforms.
\newblock {\em Mathematika}, 27(2):262--267, 1980.

\bibitem{Kec1995}
A.~S. Kechris.
\newblock {\em Classical Descriptive Set Theory}, volume 156 of {\em Graduate
  Texts in Mathematics}.
\newblock Springer-Verlag, New York, 1995.

\bibitem{Khi1997}
A.~Y. Khinchin.
\newblock {\em Continued Fractions}.
\newblock Dover Publications, Inc., Mineola, NY, {R}ussian edition, 1997.
\newblock With a preface by B. V. Gnedenko, Reprint of the 1964 translation.

\bibitem{kraaikamp2000normal}
C.~Kraaikamp and H.~Nakada.
\newblock On normal numbers for continued fractions.
\newblock {\em Ergodic Theory Dynam. Systems}, 20(5):1405--1421, 2000.

\bibitem{Lak1973}
R.~B. Lakein.
\newblock Approximation properties of some complex continued fractions.
\newblock {\em Monatsh. Math.}, 77:396--403, 1973.

\bibitem{LeV1952}
W.~J. LeVeque.
\newblock Continued fractions and approximations in {$k(i)$}. {I}, {II}.
\newblock {\em Indag. Math. (N.S.)}, 14:526--535, 536--545, 1952.
\newblock Nederl. Akad. Wetensch. Proc. Ser. A {\bf 55}.

\bibitem{Nak1976}
H.~Nakada.
\newblock On the {K}uzmin's theorem for the complex continued fractions.
\newblock {\em Keio Engineering Reports}, 29(9):93--108, 1976.

\bibitem{quas2016ergodic}
A.~Quas and T.~Soo.
\newblock Ergodic universality of some topological dynamical systems.
\newblock {\em Trans. Amer. Math. Soc.}, 368(6):4137--4170, 2016.

\bibitem{sarig1999thermodynamic}
O.~M. Sarig.
\newblock Thermodynamic formalism for countable {M}arkov shifts.
\newblock {\em Ergodic Theory Dynam. Systems}, 19(6):1565--1593, 1999.

\bibitem{AsmusSchmidt1975}
A.~L. Schmidt.
\newblock Diophantine approximation of complex numbers.
\newblock {\em Acta Math.}, 134:1--85, 1975.

\bibitem{Sch1975Subspace}
W.~M. Schmidt.
\newblock Simultaneous approximation to algebraic numbers by elements of a
  number field.
\newblock {\em Monatsh. Math.}, 79:55--66, 1975.

\bibitem{Schweiger2000}
F.~Schweiger.
\newblock {\em Multidimensional continued fractions}.
\newblock Oxford Science Publications. Oxford University Press, Oxford, 2000.

\bibitem{Walters1982}
P.~Walters.
\newblock {\em An Introduction to Ergodic Theory}, volume~79 of {\em Graduate
  Texts in Mathematics}.
\newblock Springer-Verlag, New York-Berlin, 1982.

\bibitem{Yur1995}
M.~Yuri.
\newblock Multi-dimensional maps with infinite invariant measures and countable
  state sofic shifts.
\newblock {\em Indag. Math. (N.S.)}, 6(3):355--383, 1995.

\end{thebibliography}
\end{document}